\documentclass[a4paper,11pt]{amsart}

\usepackage[utf8x]{inputenc}
\usepackage[T1]{fontenc}
\usepackage{amsfonts}
\usepackage{yfonts}
\usepackage{amsmath}
\usepackage{amsthm}
\usepackage{amssymb}
\usepackage{graphicx}
\usepackage[frenchb,english]{babel}
\usepackage{verbatim}
\usepackage{xy}
\usepackage{graphics}
\usepackage{bbm}
\usepackage{epstopdf}
\usepackage{wasysym}
\usepackage{tikz-cd}
\usetikzlibrary{arrows,decorations.pathmorphing}
\usepackage[hidelinks]{hyperref}
\usepackage[margin=3.5 cm]{geometry}

\let\~=\tilde

\newcommand{\R}{\mathbb{R}}
\newcommand{\Z}{\mathbb{Z}}
\newcommand{\N}{\mathbb{N}}
\newcommand{\Imm}{\mathrm{Im}}

\newcommand{\Tr}{\mathrm{tr}}
\newcommand{\hooklongrightarrow}{\lhook\joinrel\longrightarrow}
\newcommand{\twoheadlongrightarrow}{\relbar\joinrel\twoheadrightarrow}
\newcommand{\GT}{\mathrm{GT}}
\newcommand{\Ru}{\mathfrak{Ru}}
\newcommand{\Qu}{\mathfrak{Qu}}
\newcommand{\T}{\mathcal{T}_2(X)}
\newcommand{\Sm}{\mathcal{S}_2(X)}
\newcommand{\To}{\mathrm{T}}
\newcommand{\Le}{\mathrm{L}}
\newcommand{\m}{\mathfrak{m}}

\newtheorem{Thm}{\textbf{Theorem}}[section]
\newtheorem{Lemma}[Thm]{\textbf{Lemma}}
\newtheorem{Cor}[Thm]{\textbf{Corollary}}
\newtheorem{Prop}[Thm]{\textbf{Proposition}}

\theoremstyle{remark}

\numberwithin{equation}{section}
\newtheorem{Ex}[Thm]{\textbf{Example}}
\newtheorem{Def}[Thm]{\textbf{Definition}}
\newtheorem{Rmk}[Thm]{\textbf{Remark}}

\newtheorem{Notation}[Thm]{\textbf{Notation}}

\newtheorem{Question}{\textbf{Question}}

\newtheorem*{LFPThm}{\textbf{Lefschetz fixed point theorem}}

\title[Twisted invariants associated with bundles]{Twisted Gromov and Lefschetz invariants associated with bundles}
\date{}
\author{Gilberto Spano}

\begin{document}
\setcounter{tocdepth}{1}

\let\thefootnote\relax\footnote{Date: \today}
\let\thefootnote\relax
\let\thefootnote\relax\footnote{The author was supported by ERC LTDBud.}

\begin{abstract}
% Given a symplectic $4$-manifold $(X,\omega)$ endowed with an $\omega$-compatible almost complex structure $J$, Taubes defined a version $GT(X,\omega)$ of
% the Gromov's symplectic invariants for $(X,\omega)$, obtained by counting in a given way certain embedded $J$-holomorphic curves in $X$. A few years later
% Taubes proved that $GT(X,\omega)$ is equivalent to the Seiberg-Witten invariant of $X$, so that $GT(X,\omega)$ depends in fact only on the diffeomorphism class
% of $X$. 
Given a closed symplectic $4$-manifold $(X,\omega)$, we define a twisted version of the Gromov-Taubes invariants for $(X,\omega)$, where the twisting
coefficients are induced by the choice of a surface bundle over $X$.
%We prove that, fixed a surface bundle $F \hookrightarrow W \stackrel{\pi}{\twoheadrightarrow} X$, the
%corresponding $\pi$-twisted Gromov-Taubes series is a symplectic invariant of $(X,\omega)$.
Given a fibered $3$-manifold $Y$, we similarly construct twisted Lefschetz zeta functions associated with surface bundles: we prove that these
are essentially equivalent to the Jiang's Lefschetz zeta functions of $Y$, twisted by the representations of $\pi_1(Y)$ that are induced by
monodromy homomorphisms of surface bundles over $Y$. This leads to an
interpretation of the corresponding twisted Reidemeister torsions of $Y$ in terms of products of ``local'' commutative Reidemeister torsions.
Finally we relate the two invariants by proving that, for any fixed closed surface bundle $\mathcal{B}$ over $Y$, the corresponding twisted
Lefschetz zeta function coincides with the Gromov-Taubes invariant of $S^1 \times Y$ twisted by the bundle over $S^1 \times Y$ naturally induced by $\mathcal{B}$.

%proving that, given a closed surface bundle $F\hookrightarrow V \stackrel{\pi}{\twoheadrightarrow} Y$ over a
%fibered $3$-manifold $Y$, the corresponding twisted Lefschetz zeta function equals the twisted Gromov-Taubes invariant associated with the bundle
%$F\hookrightarrow S^1\times V \stackrel{\mathrm{id}\times\pi}{\twoheadrightarrow} S^1 \times Y$.
 
\end{abstract}

\maketitle

\tableofcontents

\section{Introduction}

To any given closed symplectic $4$-manifold $(X,\omega)$, endowed with an auxiliary generic $\omega$-compatible almost complex structure $J$,
it is possible to associate its Gromov-Taubes series $\mathrm{GT}(X,\omega,J)$ (defined in \cite{Ta1}). Even if this depends
a priori on the choice of $J$, Taubes directly showed that in fact it depends only on $X$ and the isotopy class of $\omega$. Furthermore he proved in \cite{Ta2}
that if $b_2^+(X) >1$ then $\mathrm{GT}(X,\omega)$ is equivalent to the Seiberg-Witten invariant of $X$ (see for example \cite{Morgan}), and so
it depends only on the diffeomorphism class of $X$.

Briefly, $\mathrm{GT}(X)$ can be defined in terms of a weighted count of certain embedded $J$-holomorphic surfaces $C$
in $(X,\omega,J)$. The weight of each $C$ depends on a sign and its homology class $[C]\in H_2(X,\Z)$, encoded by a formal variable
$t_{[C]}$. 

In the standard case, $\mathrm{GT}(X)$ is related to other invariants. Of particular interest for us is the case of symplectic $4$-manifolds
of the form $S^1 \times Y$: by a result of Friedl and Vidussi (\cite{F-V 1}) $Y$ must be a fibered $3$-manifold. In this case it is possible to directly
prove (see for example Section 2.6 of \cite{Hu1}) that $\mathrm{GT}(S^1 \times Y)$ coincides with the \emph{Lefschetz zeta function} $\zeta(Y)$,
which in turn is equivalent to the Reidemeister torsion of $Y$ (also called Milnor or Turaev torsion), which is a topological invariant of $Y$. 
We remark that if $L\subset S^3$ is a link, then the Reidemeister torsion of $S^3\setminus L$ is essentially equivalent to the
Alexander polynomial of $L$.

There are several refinements of the Lefschetz zeta function. One of them is the \textit{twisted Lefschetz zeta function} $\zeta_{\rho}(Y)$, defined
by Jiang (\cite{J}, \cite{J-W}), which is associated to a representation $\rho:\pi_1(Y)\rightarrow GL(n,R)$, where $R$ is a commutative ring with unit. 
Jiang proved that if $Y$ is fibered, then $\zeta_{\rho}(Y)$ is equivalent to \emph{Lin's $\rho$-twisted Reidemeister torsion of $Y$} (see for example
\cite{Lin}). Roughly speaking, if $Y$ is diffeomorphic to the mapping torus of a surface $S$ with a diffeomorphism $\phi$, $\zeta_{\rho}(Y)$
is a weighted count of periodic orbits of $\phi$, where the weight of an orbit $\delta$ depends on its Lefschetz sign, a formal variable encoding the
period of $\delta$ and the ``non-commutative'' weight $\Tr(\rho([\delta]))$, where $[\delta]$ is the conjugacy class in $\pi_1(Y)$ determined by $\delta$.

\vspace{0.3 cm}

The main motivations behind this paper are the following questions:
\begin{Question} \label{Question 1}
 There exist ``twisted versions'' of the Gromov-Taubes invariants that are related to the twisted Reidemeister torsions
 like the standard $\GT(S^1 \times Y)$ is related to the standard Reidemeister torsion of the fibered manifold $Y$?
\end{Question}
\begin{Question} \label{Question 2}
 There is a purely topological interpretation of the twisting coefficients of Jiang and Lin? Said differently, what is the topological meaning of the
 information carried by $\zeta_{\rho}(Y)$ for a given choice of $\rho$?
\end{Question}

In order to try to answer to these (apparently unrelated) questions we introduce the concept of ``bundled twistings'' for both Gromov-Taubes invariants and
Lefschetz zeta functions.

%The intuitive idea behind bundled twistings is in some sense similar to that for characteristic classes. Suppose that $\mathcal{I}$ is
%an invariant and that, given a manifold $M$, $\mathcal{I}(M)$ is defined by counting in some way a set $\{N_i\}$ of submanifolds of $M$. Given a bundle over $M$
%we can then modify $\mathcal{I}$ by taking into account a measure of how the bundle is twisted along the submanifolds $N_i$. 

We start in Section \ref{Section: Review of GT} by briefly recalling the definition of the standard Gromov-Taubes invariants.

In Section \ref{Section: Definition of the twisted invariants} we define a \emph{bundled twisted version of $\mathrm{GT}(X,\omega,J)$}, where
the aim of the twisting coefficient associated to a $J$-holomorphic curve $C \subset X$ is to detect non-commutative informations about the homotopy
class of $C$ in $X$. Roughly speaking, a bundled twisting is obtained by choosing a smooth surface bundle 
\begin{equation*} 
 F \hooklongrightarrow W \stackrel{\pi}{\twoheadlongrightarrow} X
\end{equation*}
and then twisting the weight of $C$ in $\mathrm{GT}(X,\omega,J)$ by the Gromov-Taubes invariant $\mathrm{GT}(\pi^*C)$ of
the pull-back bundle over $C$, 
%computed with respect to a natural symplectic form (see Section \ref{Section: Definition of the twisted invariants}
%for the known results about the existence of these symplectic forms)
computed with respect to formal variables encoding the homology classes in the image of
$H_2(\pi^*C,\Z) \stackrel{i_*}{\rightarrow} H_2(W,\Z)$, where $i: \pi^*C \hookrightarrow W$ is the inclusion. Since the diffeomorphism type of
$\pi^*C$ depends on the isotopy class of $C$ in $W$, it is natural to expect to get a refinement of the standard
$\mathrm{GT}(X,\omega,J)$, where only the homology class of $C$ is taken into account.

What we get is a \emph{$\pi$-twisted Gromov-Taubes series $\mathrm{GT}_{\pi}(X,\omega,J)$},
depending, a priori, on the symplectic form $\omega$ and the $\omega$-compatible almost complex structure $J$ on $X$.
In Subsection \ref{Subsection: Proof of the invariance} we will then prove the following:

\begin{Thm} \label{Theorem: GT twisted is a symplectic invariant}
 Given a smooth oriented surface bundle $(W,X,\pi,F)$ with $F$ closed, the $\pi$-twisted Gromov-Taubes series $\mathrm{GT}_{\pi}(X,\omega,J)$ is
 independent on a generic choice of $J$ and depends only on the isotopy class of $\omega$.
\end{Thm}
\begin{comment}
We observe that, while in $\mathrm{GT}(X)$ the topological information of the holomorphic curves is given only by their homology class,
in $\mathrm{GT}_{\pi}(X,\omega)$ the weights are more of homotopical nature. For instance, if two $J$-holomorphic curves $C$ and $C'$ have opposite
sign but same homology class, their total contribution to $\mathrm{GT}(X)$ is always $0$; on the other hand if $C$ is not homotopic to $C'$
it can happen that $\mathrm{GT}(\pi^*C)\neq \mathrm{GT}(\pi^*C')$ and their total contribution to $\mathrm{GT}_{\pi}(X,\omega)$ can be non zero. 
\end{comment}

%We remark that our strategy of choosing a fiber bundle to refine an invariant defined by counting holomorphic curves is original and can a priori be applied to 
%other holomorphic-curves-based invariants.

\vspace{0.3 cm}

In Section \ref{Section: Mapping tori and twisted Lesfchetz zeta functions} we apply the idea of the bundled twisting coefficients to define \emph{$\pi$-twisted Lefschetz
zeta functions} $\zeta_{\pi}(Y)$ for a fibered $3$-manifold $Y$, where now $\pi$ is the projection map of a smooth surface bundle 
\begin{equation*} \label{Equation: Diagram of the fibration in dimension 3 in introduction}
 F \hooklongrightarrow V \stackrel{\pi}{\twoheadlongrightarrow} Y.
\end{equation*}
The definition of $\zeta_{\pi}(Y)$ is conceptually similar to that of $\mathrm{GT}_{\pi}(X,\omega)$: morally, the latter is defined by counting holomorphic
curves in the pull-back bundles over $J$-holomorphic curves, while the former is defined by counting periodic orbits in the
pull-back bundles over periodic orbits. The key fact is that the pull-back bundle of $\pi$ over any periodic orbit is again a mapping torus 
with fiber $F$ and first return map determined by the \emph{monodromy homomorphism}
$$\rho_{\pi}:\pi_1(Y) \longrightarrow \mathrm{MCG}(F)$$
\emph{of the bundle}: this ``local'' mapping torus has again a standard Lefschetz zeta function, that can be used as twisting weight for the underlying
orbit in $Y$.

It turns out that our $\pi$-twisted Lefschetz zeta function is essentially equivalent to the Jiang's Lefschetz zeta function, ``algebraically''
twisted with respect to the matrix representations
$${\rho_{\pi}}_*:\pi_1(Y) \longrightarrow \mathrm{GL}\left(H_*(F)\right)$$
induced in homology by $\rho_{\pi}$
(see Theorem \ref{Theorem: Relation between geometric and algebraic twistings for the zeta function.} below). As a corollary we obtain that for any 
surface bundle $(V,Y,\pi,F)$, $\zeta_{\pi}(Y)$ is a \emph{topological invariant}.
In fact, since $\zeta_{\pi}(Y)$ is defined by counting orbits (in the total space $V$) just with sign and their homological weight, it can be seen
as a purely geometrical interpretation of the corresponding twisted Reidemeister torsion of Lin. The intuitive idea behind this phenomenon
is the following: using Lefschetz fixed point theorem, we can interpret the Jiang's twisting coefficients $\Tr({\rho_{\pi}}_*([\delta^n]))$ of
the iterates of an orbit $\delta$ as a signed sum of periodic points (orbits) in the mapping torus $\pi^*\delta$ of $(F,\rho_{\pi}([\delta]))$,
which in turn gives (a version of) the Lefschetz zeta function of $\pi^*\delta$. This provides an answer to Question \ref{Question 2} for
the family of \emph{bundled representations} of $\pi_1(Y)$ of Definition \ref{Definition: Bundled representations} below (cf. Example
\ref{Example: All representations on GL 2g Z are bundled} and Remark \ref{Remark: What are the bundled representations?}).

\begin{comment}
\begin{Thm} \label{Theorem: Relation between geometric and algebraic twistings for the zeta function in intro}
 Let $(V,Y,\pi,F)$ be a surface bundle like in Diagram (\ref{Equation: Diagram of the fibration in dimension 3 in introduction}), with $Y$
 fibered $3$-manifold. For, $i = 0,1,2$, let 
 $$\m_i = {(\m_{\pi}})_i: \pi_1(Y) \rightarrow GL(H_i(F,\Z))$$
 be the group representations induced in homology by the monodromy map of the bundle. Then:
 \begin{equation}
  \zeta_{\pi}(Y) = \prod_{i=0}^2 (\zeta_{\m_i}(Y))^{(-1)^i}.
 \end{equation}
\end{Thm}
\end{comment}

\vspace{0.3 cm}

The reason for which the two motivational questions above are related is explained by the following:
\begin{Thm} \label{Theorem: GT twisted = Zeta twisted in intro}
 Let $Y$ be the mapping torus of a closed surface $S$ and a diffeomorphism $\phi$ and let $(V,Y,\pi,F)$ be a smooth oriented surface bundle with
 $F$ closed. Let
 $$(S^1 \times V,S^1 \times Y,\mathrm{Id} \times \pi,F)$$
 be the natural bundle induced by the product by $S^1$ and consider the symplectic form $\omega = \Omega + ds \wedge dt$ on $S^1 \times Y$, where $\Omega$ is
 any symplectic form on $S$ and $t$ and $s$ are coordinates for $[0,1]$ and, respectively, $S^1$. Then:
 \begin{equation*}
  \mathrm{GT}_{\mathrm{Id} \times \pi}(S^1 \times Y,\omega) = \zeta_{\pi}(Y).
 \end{equation*}
\end{Thm}

\vspace{0.3 cm}

By last theorem, and the relation between $\zeta_{\pi}$ and Lin's twisted Reidemeister torsion, $\mathrm{GT}_{\pi}(X,\omega)$ provides
an answer to Question \ref{Question 1} and can be considered as a ``non-triviality'' result for the twisted Gromov-Taubes invariants. Moreover,
since twisted Reidemeister torsions are topological invariants, it is natural to ask whether the twisted $\mathrm{GT}$'s are in fact
in general independent also on the choice of the symplectic form and how powerful they are in detecting generic closed symplectic $4$-manifolds. 

\vspace{0.3 cm}

We remark that twisted Reidemeister torsions are in general much stronger invariants than the standard version. For example it has been proved that
they detect fibered (non-necessarily closed) $3$-manifolds (\cite{F-V 1}), Thurston norms and genus of knots (\cite{F-V 3}) and, very recently, 
the hyperbolic volume of knot exteriors (\cite{Goda}). Also, they give useful tools in detecting sliceness and knot concordance: see
\cite{F-V 3} for these and other applications. It is then natural to expect to find several useful applications also for the twisted Gromov-Taubes invariants.

%\addtocontents{toc}{\protect\setcounter{tocdepth}{1}}
\subsubsection*{Acknowledgements} The author is grateful to Paolo Ghiggini and András Stipsicz for many helpful suggestions and discussions.
He also would like to thank Antonio Alfieri, Vincent Colin and Chris Wendl for the interest they showed about this work.
%\addtocontents{toc}{\protect\setcounter{tocdepth}{2}}

\section{Review of the Gromov-Taubes invariants} \label{Section: Review of GT}

In this brief review of $\mathrm{GT}(X,\omega)$ we will follow \cite[Chapter 2]{Hu1}, \cite{I-P} and \cite{Ta1}.

Fix $(X,\omega)$ and an $\omega$-compatible almost complex structure $J$ over $X$.
A \emph{$J$-holomorphic curve} in $(X,J)$ is a smooth map $u : (\Sigma,j) \rightarrow (X,J)$ from a connected compact Riemann surface $(\Sigma,j)$ such that 
$$\overline{\partial}f := \frac{1}{2}(du \circ j - J \circ du) = 0$$
and considered up to the following equivalence relation: given $u_i : (\Sigma_i,j_i) \rightarrow (X,J)$, $i =1,2$, then $u_1$ and $u_2$ are equivalent
if there exists a biholomorphism $\phi :(\Sigma_1,j_1) \rightarrow (\Sigma_2,j_2)$ such that $u_2 \circ \phi = u_1$. We will usually denote the image of $u$ by
$C_u$.

If the holomorphic curve is connected it can be of two types:
\begin{itemize}
 \item \emph{somewhere injective}, i.e. $\exists x \in \Sigma$ such that $u^{-1}(u(x)) = \{x\}$ and $d_x u$ is injective; 
 \item \emph{multiply covered}, i.e. $u$ factors through a degree bigger than $1$ branched cover over another holomorphic curve. 
\end{itemize}
One can show that in the first case $C_u$ is an embedded surface outside a finite number of points, called \textit{singularities}.
%A singularity is called a
%\textit{node} of $u$ if it sits on a transverse self-intersection of $C_u$ with total multiplicity $2$. In the rest of the paper we will always
%assume that every singular point is a node.

Holomorphic curves in dimension $4$ have some nice properties. One of them is the \textit{positivity of intersection}, which can be formulated as
follows. Let $u$ and $v$ be two distinct holomorphic curves in $(X,J)$. %Then $\#(C_u \cap C_v) < \infty$. Moreover,
If $P \in C_u \cap C_v$, then its contribution $m_P$ to the algebraic intersection number
$ C_u \cdot C_v$ is strictly positive, and $m_P = 1$ if and only
if $u$ and $v$ are embeddings near $P$ that intersect transversely in $P$.

Note that the assumption that $u_1$ and $u_2$ are distinct is crucial: for example one can have embedded holomorphic spheres with negative
autointersection. 

Another relevant property is the following \emph{adjunction formula}. Let $u : (\Sigma,j) \rightarrow (X,J)$ be a somewhere injective $J$-holomorphic
curve. Then:

\begin{equation}
 \langle c_1(TX),[C_u]\rangle = \chi(\Sigma) + [C_u] \cdot [C_u] -2\delta(C_u)
\end{equation}

\phantom{}
\noindent where $\delta$ is a weighted count of the singularities of $C_u$. In particular $C_u$ must have genus
\begin{equation}
 g_{C_u} = 1+\frac{1}{2}\big([C_u] \cdot [C_u] - \langle c_1(TX),[C_u]\rangle -\delta(C_u)\big).
\end{equation}
%which, of course, coincides also with the genus of $\Sigma$.

\begin{Def}
 Given $A \in H_2(X,\Z)$, we define $g_A$ to be the genus of any embedded $J$-holomorphic curve in the homology class $A$ as given by the adjunction
 formula:
 \begin{equation*}
  g_A := 1+\frac{1}{2}\big(A \cdot A - \langle c_1(TX),A\rangle\big).
 \end{equation*}
\end{Def}

In the definition of the Gromov-Taubes invariants a special role will be played by the homology classes $A$ with
$g_A=0$ and $g_A=1$. We define then the following subsets of $H_2(X,\Z)$:
\begin{equation*}
 \begin{array}{lll}
  \Sm & := & \{A \in H_2(X,\Z)\ |\ g_A =0,\ A \cdot A = 0 \};\\
  \T & := & \{A \in H_2(X,\Z)\ |\ g_A =1,\ A \cdot A = 0 \}.
 \end{array}
\end{equation*}

To any $A \in H_2(X,\Z)$ we can also associate the following crucial number:
\begin{equation*}
 d_A := A \cdot A + \langle c_1(TX),A\rangle .
\end{equation*}
\noindent This (even) number coincides with the expected dimension of the moduli space $\mathcal{M}_A(X,J)$ of embedded holomorphic curves
in the homology class $A$, so that, whenever $d_A = 0$, we are able to count the elements of $\mathcal{M}_A(X,J)$. In general, if $d_A \geq 0$ it
is possible to cut the dimension of $\mathcal{M}_A(X,J)$ by restricting the attention to the $0$-dimensional subspace
$\mathcal{M}_{A}(X,J,\Omega_{d_A})$ of $\mathcal{M}_A(X,J)$ consisting of those holomorphic curves that pass through a fixed set $\Omega_{d_A}$ of
$d_A/2$ generic points in $X$. Obviously if ${d_A}=0$ then $\mathcal{M}_{A}(X,J,\Omega_0) =\mathcal{M}_{A}(X,J,\emptyset) = \mathcal{M}_A(X,J)$.

\subsection{Definition of $\GT$} \label{Subsection: Definition of GT} $ $\\
\noindent The Gromov-Taubes invariant can be defined as a weighted count of the elements of the moduli spaces $\mathcal{M}_{A}(X,J,\Omega_{d_A})$ with $A\in H_2(X,\Z)$
such that $d_A \geq 0$; the weight of each embedded $J$-holomorphic curve is defined in terms of: 
\begin{enumerate}
 \item the sign $\epsilon(C_u)$ of the determinant of a differential operator associated to each $C_u$ (see \cite{Ru} or \cite{Ta1} for more
  details);
 \item formal variables $t_A,A\in H_2(X)$ keeping track of the homology class of the curves and satisfying $t_{A+B} = t_A\cdot t_B$.
\end{enumerate}
The invariant can then be seen as a series
$$\mathrm{GT}(X) = \sum_A \mathrm{GT}(X,A) t_A \in \Z[[ H_2(X)]].$$

\begin{Def} \label{Definition: Decompositions D of A}
 For any $A \in H_2(X)$, let $\mathcal{D}(A)$ be the set of couples $(A_i,n_i) \in H_2(X) \times \N^*$ such that:
\begin{enumerate}
 \item $A = \sum_i n_i A_i$;
 \item $d_{A_i} \geq 0$;
 \item the $A_i$'s are distinct and not \textit{multiply toroidal} (i.e. of the form $m B$ with $m \geq 2$ and $g_B =1$);
 \item $A_i \cdot A_j = 0$ when $i\neq j$;
 \item if $A_i\cdot A_i \neq 0$ then $n_i =1$;
 \item if $A_i\cdot A_i =0$ then $n_i$ can be any positive integer.
\end{enumerate}
 \end{Def}
Observe that if there exists a holomorphic curve in a class $A_i$ with $d_{A_i} \geq 0$ and $A_i\cdot A_i =0$, then the adjunction formula implies
that $A_i$ is either in $\Sm$ or in $\T$. Note moreover that
$t_A =\sum_i t_{n_i A_i} = \sum_i t_{A_i}^{n_i}$ and $d_A=\sum_i n_i d_{A_i}$ (cf. \cite[Lemma 2.4]{Hu1}).

Now, given $A$, fix $d_A/2$ generic points in $X$; $\mathrm{GT}(X,A)$ is then of the form:
\begin{equation*} %\label{Equation: GT(X,A)}
 \mathrm{GT}(X,A) = \sum_{(A_i,n_i) \in \mathcal{D}(A)} d_A! \left(\prod_{A_i \notin \T} \frac{\mathfrak{Ru}(A_i)^{n_i}}{d_{A_i}!n_i!} \cdot
 \prod_{A_i \in \T} \mathfrak{Qu}(A_i,n_i)\right).
\end{equation*}

Note that in the first product if $n_i \neq 1$ then $A_i \in \Sm$. For any $B \in H_2(X)$, the quantity $\mathfrak{Ru}(B)$ is the
\textit{Ruan invariant of $(X,B)$} (\cite{Ru}) and it is defined by:

\begin{equation}\label{Equation: Ruan invariant}
 \mathfrak{Ru}(B) := \sum_{C_u \in {\mathcal{M}}_{B}(X,J,\Omega_{d_B})} \epsilon(C_u). 
\end{equation}

It is important here to recall that each of the $C_u$ is embedded and connected by definition.

The definition of $\mathfrak{Qu}(B,n)$ for a (non multiply) toroidal class $B \in \T$ is a bit more complicated. It can be expressed in
the form

\begin{equation} \label{Equation: Qu(B,m)}
 \mathfrak{Qu}(B,m) := \sum_{\{(C_k,m_k)\}}\prod_k r(C_k,m_k)
\end{equation}

\noindent where the sum is over the sets of couples $(C_k,m_k) \in \mathcal{M}_{c_k B}(X,J) \times \N^*$ with $c_k\geq 1$ and
$\sum_k c_k m_k = m$.
Observe that $B \in \T$ if and only if $ B \cdot B = \langle c_1(TX),B\rangle =0$, so that also $cB\in \T$
for any $c \geq 1$. Then all the $C_k$ are (embedded) tori.

If $C$ is an embedded holomorphic torus $C$ with homology class in $\T$ , in \cite{Ta1} Taubes extends the definition of the sign $\epsilon$ also
to the holomorphic double covers of $C$. The numbers $r(C,l)$ can be defined in terms of the \emph{Taubes' generating functions} $P(C,z)$: these
depend on the signs of $C$ and of its three connected double covers, which can be identified with the cohomology classes
$\iota_1,\iota_2,\iota_3 \in H^1(C,\Z/2)$.

%The numbers $r(C,l)$ can be defined explicitly as in \cite[Definition 3.2]{Ta1}, but also via the \emph{Taubes' generating functions}. 
We say that $C$ is of \textit{type} $(\epsilon,s) \in \{\pm\} \times \{0,1,2,3\}$ if $\epsilon(C) = \epsilon$ and exactly $s$ of the
connected double covers $\iota_1,\iota_2,\iota_3$ have sign $-$. We define then
\begin{equation*}
 P(C,z) := P_{\epsilon(C),s}(z).
\end{equation*}
where the functions $P_{\pm,s}(z)$ are defined in terms of $P_{+,0}(z)$ by

\begin{equation} \label{Equation: Generating functions standard}
 \begin{array}{lcl}
  P_{+,1}(z) := \dfrac{P_{+,0}(z)}{P_{+,0}(z^2)} & , & P_{+,2}(z) := \dfrac{P_{+,0}(z)P_{+,0}(z^4)}{(P_{+,0}(z^2))^2}\ ,\\ 
   & & \\
  P_{+,3}(z) := \dfrac{P_{+,0}(z)P_{+,0}(z^4)}{(P_{+,0}(z^2))^3} & , & P_{-,s}(z) := \dfrac{1}{P_{+,s}(z)}
 \end{array}
\end{equation}
\phantom{}\\
\noindent and then normalized by setting
\begin{equation} \label{Equation: Taubes normalization of P_+0}
 P_{+,0}(z) := \dfrac{1}{1-z}.
\end{equation}
\noindent In Section \ref{Subsection: Proof of the invariance} we will say something about where these generating functions come from.
%Intuitively, the formal variable
%$z$ will keep track of the homology class of the embedded tori and its exponent of the multiplicity with which the tori are counted.

For any $l$, the number $r(C,l)$ is then defined as the coefficient of $z^l$ in the formal power series expansion (about $0$) of $P(C,z)$. Then
$\mathfrak{Qu}(B,m)$ in Equation (\ref{Equation: Qu(B,m)}) is the coefficient of $t_{mB}$ in
\begin{equation*}
 \prod_{C \mbox{\scriptsize{ embedded }}|\ [C] \in \T} P(C,t_{[C]}).
\end{equation*}

\begin{Rmk} \label{Remark: Meaning of the terms in the powers of z}
 Observe that in the definition of Taubes' generating functions there are terms that depend on $z$, $z^2$ and $z^4$: intuitively the term in $z^i$
 counts (with signs) $i$-fold covers of the tori. Note moreover that using these functions, for any torus $C$ we get $r(C,1) = \epsilon(C)$ as expected. 
\end{Rmk}

\subsection{About the symplectic invariance of $\GT$} \label{Subsection: About the symplectic invariance of GT} $ $ \\
\indent The definition of $\GT$ depends strongly on the choice of the symplectic form $\omega$, the $\omega$-compatible almost complex structure and
the sets $\Omega_d$.
In Sections 4 and 5 of \cite{Ta1}, Taubes proves that, for a fixed $\omega$, there exists a dense open subset $\mathcal{U}_d$ of the
set of the ``\textit{admissible}'' couples $\mathcal{A}_d \subset \{(J,\Omega_d)\}$
such that $\mathcal{M}_A(X,J,\Omega_d)$ (and so also both $\Ru$ and $\Qu$) is independent on $(J,\Omega_d)$ whenever
$(J,\Omega_d)$ varies within a small enough neighborhood of any point of $\mathcal{U}_d$. 
Moreover (see assertion 3 of \cite[Proposition 5.2]{Ta1}):
\begin{Lemma} \label{Lemma: A path of symplectic structures gives a smooth fibered manifold}
 For any smooth path $\{\omega_t\ |\ t\in [0,1]\}$ of symplectic structures on $X$, there exists a smooth path $\{J_t\ |\ t\in [0,1]\}$ of
 $\omega_t$-compatible almost complex structures (with prescribed $J_0$ and $J_1$) such that the fibered product
 \begin{equation} \label{Equation: Simple 1-dimensional space of curves given by a path of sympl. str.}
 \mathcal{C}_A(X,\{J_t\},\Omega_d) := \{(t,C_t)\ |\ t\in [0,1],\ C_t \in \mathcal{M}_A(X,J_t,\Omega_d)\} 
\end{equation}
has the structure of $1$-dimensional oriented manifold, which is also compact if $A$ is not multiply toroidal.
\end{Lemma}
This is enough to conclude that $\Ru$ is a symplectic invariant. 

\vspace{0.3 cm}

To prove that also $\Qu$ is a symplectic invariant, Taubes needed to menage the case of multiply toroidal classes. The issue with these classes is
compactness: if $A \in \T$ is not multiply toroidal and $m$ is a positive integer, then it can happen that a sequence of embedded tori in
$\mathcal{C}_{mA}(X,\{J_t\})$ has no limit in $\mathcal{C}_{mA}(X,\{J_t\})$ (we avoid to refer to $\Omega_0 = \emptyset$ in the notation).
To solve the problem of the sequences that have no ``true'' limits, for any fixed $n$, Tubes considers the set
\begin{equation} \label{Equation: Union K of 1-dimensional spaces of tori given by a path of sympl. str.}
 \mathcal{K}_{nA}(X,\{J_t\}) := \bigcup_{m \leq n} \mathcal{C}_{mA}(X,\{J_t\})
\end{equation}
endowed with the topology induced by the disjoint union. Note that, by definition, every point in $\mathcal{K}_{nA}(X,\{J_t\})$ corresponds
to an \emph{embedded torus} whose homology class is some multiple of $A$.
Taubes proves then the following Lemma (\cite[Lemma 5.8]{Ta1}). 

\begin{Lemma} \label{Lemma: Union of 1-dimensional spaces of tori given by a path of sympl. str. is a manifold}
 Even if a sequence $\{(t,C_t)\}_{t\rightarrow t_0}$ in $\mathcal{C}_{mA}(X,\{J_t\})$ has no limit in $\mathcal{C}_{mA}(X,\{J_t\})$,
 $(t,C_t)$ still has a ``\emph{weak limit}'' in $\mathcal{K}_{nA}(X,\{J_t\})$ for any $n\geq m$, i.e. there exist $p$ and $q$ positive integers
 with $pq=m$ such that $\{C_t\}_{t\rightarrow t_0}$ converges to a $p$-fold holomorphic cover of a torus belonging to
 $\mathcal{M}_{qA}(X,J_{t_0})\subset\mathcal{C}_{qA}(X,\{J_t\})\subset \mathcal{K}_{nA}(X,\{J_t\})$.
\end{Lemma}

If the path $\{\omega_t\ |\ t\in [0,1]\}$ is ``reasonable'', Taubes proves that we can have only weak limits with $p=2$ and that the
total number of the corresponding ``bifurcation points'' $(t_0,C_{t_0})$ in $\mathcal{K}_{nA}(X,\{J_t\})$ is finite: here is where
the double covers of the tori come into play (see Lemmas 5.8-5.11 in \cite{Ta1}). Taubes studies what happens at each bifurcation point
to the signs of the tori and of their double covers and finds that, in order to produce a symplectic invariant by counting tori $C$, their weights $P(C)$
must satisfy certain relations induced by the path $\{\omega_t\}$ when it crosses bifurcation points. These relations are exactly those in
\ref{Equation: Generating functions standard}. The normalization (\ref{Equation: Taubes normalization of P_+0}) for $P_{+,0}$ has then be chosen by
Taubes to make $\GT$ coincide with the Seiberg-Witten invariants.

\section{Twisted Gromov-Taubes invariants} \label{Section: Definition of the twisted invariants}

Fix a closed symplectic $4$-manifold $(X,\omega)$ endowed with an $\omega$-compatible almost complex structure $J$ and fix a smooth fiber bundle
\begin{equation} \label{Equation: Diagram of the fibration in dimension 4}
 F \hooklongrightarrow W \stackrel{\pi}{\twoheadlongrightarrow} X
\end{equation}
where $F$ is a closed oriented surface. Our aim is to define an analogue of $\GT(X,\omega)$ in which the weight of each $J$-holomorphic embedded curve
$C \subset X$ is (morally) twisted by the Gromov-Taubes invariant $\GT(\pi^*C)$ of the total space $\pi^*C$ of the restriction of $\pi$ to $C$ (here we
prefer to use the notation $\pi^*C$ instead of $\pi^{-1}(C)$). 

\begin{Thm}[Thurston, \cite{Th1}]\label{Theorem: Thurston} Given a smooth surface bundle
 \begin{equation*} 
  F \hooklongrightarrow N \stackrel{\pi}{\twoheadlongrightarrow} M
 \end{equation*}
 over a $2n$-dimensional symplectic manifold $(M,\omega)$, if the homology class of the fiber in $H_2(N,\R)$ is non-zero, then there exists a closed
 $2$-form $\alpha$ on $N$ that is non-singular on each fiber and such that 
 \begin{equation}\label{Equation: Thurston's symplectic form}
  \alpha + \pi^*\omega \in \Omega^2(N)
 \end{equation}
 is symplectic. 
\end{Thm}

Remark that the $2$-form in (\ref{Equation: Thurston's symplectic form}) makes all the fibers symplectic and that the condition 
$[F] \neq 0$ in $H_2(N,\R)$ is always satisfied when the genus of $F$ is greater than $1$. 
When the dimension of the base $M$ is $2$ we have also the following:
\begin{Lemma} \label{Lemma: The acs of the Thurston symplectic forma are the same}
 Assume that in the last theorem $n=1$, so that $(M,\omega)$ is a symplectic surface. Let $\omega_N$ and $\omega_N'$ be two symplectic forms on $N$ like in 
 (\ref{Equation: Thurston's symplectic form}). Then the homotopy classes of the $\omega_N$- and $\omega_N'$-compatible almost complex structures
 coincide.
\end{Lemma}
\proof
 Since the fibers are symplectic with respect to $\omega_N$ and $\omega_N'$ and any two positive symplectic forms on an oriented surface are
 deformation
 equivalent, we can find $\omega_N$- and, respectively, $\omega_N'$-compatible almost complex structures $J_N$ and $J_N'$ whose restriction to
 the tangent space of the fiber $T_x F \subset T_x N$ coincide for any $x \in N$. Let $O_x$ and $O_x'$ be the orthogonal
 complements of $T_x F$ in $T_x N$ with respect to the Riemannian metrics $\omega_N(\cdot,J\cdot)$ and, respectively, $\omega'_N(\cdot,J'\cdot)$.
 These give two distributions $O$ and $O'$ of ($J$- and, respectively, $J'$-holomorphic) tangent planes in $T N$, which are homotopic
 since the space of Riemannian metrics over $N$ is contractible. This implies that there is a homotopy between the two decompositions
 $TF \oplus O$ and $TF \oplus O'$ of $TN$, which sends $J$ to $J'$.
\endproof

In order to define the $\pi$-twisted Gromov-Taubes invariants, we will define first $\pi$-twisted versions of $\Ru$ and $\Qu$. For simplicity,
from now on and if not stated otherwise, all the homology groups will be considered with coefficients in $\Z$. Moreover we will often
implicitly assume that our surface bundles have non-trivial homology class of the fiber.

\subsection{Twisted $\Ru$} \label{Subsection: Twisted Ru} 
\begin{Notation} \label{Notation: Homology classes B_W}
 If $M$ is a manifold, $L\subset M$ a submanifold and $i\in \N$, given a homology class $A \in H_i(L)$ we will call
 $A_M$ its image in $H_i(M)$ under the homomorphism induced by the inclusion $L \hookrightarrow M$. 
 Moreover, given a group $(G,+)$ and a ring $R$, $R[[G]]$ will denote the polynomial ring $R[[\{t_g|g \in G\}]]$ with the usual relations
 $t_g \cdot t_{g'} = t_{g + g'}$. When $G$ will be given as a direct sum $G = H \oplus I$, in order to keep distinct the variables associated to
 elements in $H$ and $I$, $\Z[[G]]= \Z[[H \oplus I]]$ should be thought of as $R[[\{t_h|h \in H\}]][[\{t_i|i \in I\}]]$.
\end{Notation}

For the rest of this section we fix a closed symplectic $4$-manifold $(X,\omega)$, an $\omega$-compatible almost complex structure
$J$ and a smooth surface bundle $(W,X,\pi,F)$ with $[F] \neq 0 \in H_2(W,\R)$.

Given a surface $C \subset X$, consider the $4$-manifold $\pi^*C$. The natural inclusion $\pi^*C \hookrightarrow W$ induces in homology the 
homomorphism $B \mapsto B_W$, where $B \in H_2(\pi^*C)$. Considering formal variables $z_D$, for $D \in H_2(W)$, analogue to the $t_A$'s of last
section, we have then the ring homomorphism
\begin{equation*}
 \begin{array}{ccc}
  \Z[[H_2(\pi^*C)]] & \longrightarrow & \Z[[H_2(W)]].\\
      at_B          & \longmapsto     & az_{B_W}
 \end{array}
\end{equation*}

Now, if $C \subset X$ is a $J$-holomorphic embedded surface, $\omega|_{C}$ is symplectic and by Theorem \ref{Theorem: Thurston} there exists
a closed $\alpha_C \in \Omega^2(\pi^*C)$ such that the $2$-form 
\begin{equation} \label{Equation: Symplectic form above the curves}
 \omega_C =\alpha_C + \pi^*\omega|_{C} \in \Omega^2(\pi^*C)
\end{equation}
is symplectic and makes the fibers symplectic as well. We can then consider the Gromov-Taubes invariants
\begin{equation}
 \mathrm{GT}(\pi^*C,\omega_C) = \sum_{B \in H_2(\pi^*C)} \mathrm{GT}(\pi^*C,\omega_C,B) t_{B} \in \Z[[H_2(\pi^*C)]].
\end{equation}
By the statement of the equivalence between Gromov-Taubes and Seiberg-Witten invariants, $\mathrm{GT}(\pi^*C,\omega_C)$ depends on $\omega_C$
only through the $\mathrm{Spin}^c$-structure of an $\omega_C$-compatible almost complex structure (see \cite{Ta2}). Lemma
\ref{Lemma: The acs of the Thurston symplectic forma are the same}, implies then that $\mathrm{GT}(\pi^*C,\omega_C)$ does not depend on the particular
$\omega_C$ given by Theorem \ref{Theorem: Thurston}.

%Fixed a symplectic structure $\omega_W$ on $W$ like in (\ref{Equation: Thurston's symplectic form}), let $\omega_C :=i^*\omega_W$ be the symplectic form on
%$\pi^*C$ given by Lemma \ref{Lemma: pi*C is canonically symplectic} and consider the series
%\begin{equation*}
% \sum_{B \in H_2(\pi^*C)} \mathrm{GT}(\pi^*C,\omega_C,B) z_{B_W} \in \Z[[H_2(W)]].
%\end{equation*}
%This series depends on $\omega_W$; in fact the equivalence between Seiberg-Witten and Gromov-Taubes invariants gives
%\begin{equation} \label{Equation: SW=GT for pi^*C}
% \mathrm{GT}(\pi^*C,\omega_C,B) = \mathrm{SW}(\pi^*C,\mathfrak{s}_{\omega_C} + \mathrm{PD}(B))
%\end{equation}
%where $\mathfrak{s}_{\omega_C}$ is the canonical $\mathrm{Spin}^c$-structure induced by $\omega_C$ and $\mathrm{PD}$ denotes the Poincaré dual.
%If $\omega'_W$ is another symplectic form on $W$ and $\omega'_C := i^*\omega'_W$, then for any $B \in H_2(\pi^*C)$, the relation
%(\ref{Equation: SW=GT for pi^*C}) gives: 
%\begin{equation} \label{Equation: GT pi^*C changing spin^c structure}
% \mathrm{GT}(\pi^*C,\omega'_C,B) = \mathrm{GT}(\pi^*C,\omega_C,\mathrm{PD}(\mathfrak{s}_{\omega'_C} - \mathfrak{s}_{\omega_C})+B).
%\end{equation}

%We can normalize $\mathrm{GT}(\pi^*C,\omega_C,B)$ by fixing a $\mathrm{Spin}^c$-structure over $\pi^*C$. The natural choice is the restriction
%$\pi^*\mathfrak{s}_{\omega}$ to $\pi^*C$ of the pull-back via $\pi$ of the $\mathrm{Spin}^c$-structure $\mathfrak{s}_{\omega}$ over $X$
%determined by $\omega$.

\begin{Def} \label{Definition: GT_W for embedded surfaces}
 Given $(W,X,\pi,F)$ and a $J$-holomorphic embedded surface $C \subset X$ we define
 %If $\omega_W$ is a
 %symplectic structure in $W$ obtained via Theorem \ref{Theorem: Thurston} and $\omega_C$ is the induced symplectic structure on $\pi^*C$, we define:
% $$\GT_W(\pi^*C,B) := \displaystyle  \mathrm{GT}(\pi^*C,\omega_C,\mathrm{PD}(\pi^*\mathfrak{s}_{\omega} - \mathfrak{s}_{\omega_C})+ B)$$
% \noindent and
  $$\GT_W(\pi^*C) := \displaystyle \sum_{B \in H_2(\pi^*C)} \mathrm{GT}(\pi^*C,\omega_C,B) z_{B_W}\ \in\ \Z[[H_2(W)]],$$
  where $\omega_C$ is any symplectic form like in (\ref{Equation: Symplectic form above the curves}).
\end{Def}

%where,
%The notation above is used to remind that $\zeta_{\mathrm{A}}(\phi)$ should be thought of as a function in the formal variable $t$ and
%coefficients in $\Z[[H_1\left(\To_{\phi}\right)]]$. 

Observe that $\GT_W(\pi^*C)$ depends only on the isotopy class of $C$ in $X$. Now, for any
$d\in \N$, fix a set $\Omega_{d}$ of $d/2$ generic points in $X$ and, for a homology class $A$, let
$\mathcal{M}_{A}(X,J,\Omega_{d_A})$ be as in the definition of the standard $\mathfrak{Ru}(A,\omega,J,\Omega_{d_A})$.

\begin{Def} \label{Definition: Ruan invariant twisted}
 Given the bundle $(W,X,\pi,F)$ and a class $A \in \left(H_2(X) \setminus \T\right)$, we define the \emph{$\pi$-twisted Ruan invariant} of $(X,A)$ by:
 \begin{equation}
  \mathfrak{Ru}_{\pi}(A,\omega,J,\Omega_{d_A}) := \sum_{C \in \mathcal{M}_{A}(X,J,\Omega_{d_A})} \epsilon(C_u) \cdot \GT_W(\pi^*C).
 \end{equation}
\end{Def}
\noindent When $\omega$, $J$ and $\Omega_d$ are understood we will simply write $\Ru_{\pi}(A)$ for $\mathfrak{Ru}_{\pi}(A,\omega,J,\Omega_d)$.

%\begin{Rmk}
% Observe that for another $\omega'_W$ obtained from $\omega$ via Theorem \ref{Theorem: Thurston}, $\Ru_{\pi}(A,\omega_W',J,\Omega_d)$ differs from
% $\Ru_{\pi}(A,\omega_W,J,\Omega_d)$ by an overall shift, since Equation (\ref{Equation: GT pi^*C changing spin^c structure}) implies that
% $$\Ru_{\pi}(A,\omega_W',J,\Omega_d) = \Ru_{\pi}(A,\omega_W,J,\Omega_d) .$$
% In fact we could actually replace the definition at point (2) by a similarly defined Seiberg-Witten version
% $\mathrm{SW}_W(\pi^*C)$ and in what follows we would still get a symplectic invariant.
 %On the other hand, here we want to define an invariant via a pure count of holomorphic curves. 
 %Moreover, because of Theorem \ref{Theorem: Classification of symplectic surface bundles}, ``almost
 %all'' $\pi^*C$ admit a symplectic structure.
 
% It is important to keep in mind that the only
% symplectic information that interests us here is about $X$: the aim of our twisting coefficients induced by the bundle is to catch pure
% topological informations about the holomorphic curves in $X$. 
%\end{Rmk}

\subsection{Twisted $\Qu$} \label{Subsection: Twisted Qu}  $ $ \\
\indent As in the standard case, the count for homology classes in $\T$ is more complicated. As recalled in last section, the weight in $\GT(X)$
for a $J$-holomorphic embedded torus $C$ with $[C]\in \T$ depends on the signs $\epsilon(C_{\iota_i})$ of its connected double covers
$C_{\iota_i}$, $i \in \{1,2,3\}$.

We want to define $\pi$-twisted analogues of the Taubes' generating functions in (\ref{Equation: Generating functions standard}).
To do that
%we need to understand in more details the functions $P_{\pm,i}(z)$. 
it will be convenient to have an explicit description of the double covers $\iota_i$ of $C$. Identify $C$ with $\R^2/\Z^2$ and fix the basis of
$H_1(C,\Z)$ induced by the ordered couple of segments $([0,1]\times \{0\},\{0\} \times [0,1])$ of $\R^2$ quotiented by $\Z^2$.
The covering spaces $C_{\iota_i}$, $i \in \{1,2,3\}$, of the three relevant double covers can then be represented as the images on $\R^2$ of
the square $[0,1]^2$ under the linear maps

\begin{eqnarray*} %\label{Equation: The maps induced ba the iotas in homology}
 {l_{\iota_1}} := \left(\begin{array}{cc}
                         2 & 0 \\
                         0 & 1
                         \end{array}\right)\ ,\ 
                         &                     
 {l_{\iota_2}} := \left(\begin{array}{cc}
                         1 & 0 \\
                         0 & 2
                         \end{array}\right)\ ,\ 
                         &
 {l_{\iota_3}} := \left(\begin{array}{cc}
                         1 & -1 \\
                         1 & 1
                         \end{array}\right),
\end{eqnarray*}
with their boundary quotiented by the relation induced by the quotient $[0,1]^2/\Z^2$. Observe that the covering map
$f_{\iota_i}:C_{\iota_i} {\longrightarrow} C$ is given by the restriction to $\Imm(l_{\iota_i})$ of quotient of $\R^2$ by $\Z^2$.
The relevant covers can then be identified with the elements $\iota_1,\iota_2,\iota_3 \in H^1(C,\Z/2) \cong \mathrm{Hom}\left(H_1(C,\Z),\Z/2\right)$
represented by the matrices
\begin{eqnarray} \label{Equation: The iotas in terms of cohomology classes}
 \iota_1:=\begin{pmatrix}
           1 &  0
          \end{pmatrix}\ ,\ 
 \iota_2:=\begin{pmatrix}
           0 &  1
          \end{pmatrix}\ ,\ 
 \iota_3:=\begin{pmatrix}
           1 &  1
          \end{pmatrix},
\end{eqnarray}
acting on the left on column vectors with respect to the fixed basis of $H_1(C,\Z)$.

Note that the image under $l_{\iota_i}$ of the couple of segments $([0,1]\times \{0\},\{0\} \times [0,1])$ induces a basis of $H_1(C_{\iota_i},\Z)$
and that the covering map $f_{\iota_i}$ induces in homology the map
${f_{\iota_i}}_* : H_1(C_{\iota_i},\Z) \longrightarrow H_1(C,\Z)$, which, with respect to the fixed bases, is given exactly by the matrices above:

\begin{eqnarray} \label{Equation: The maps induced ba the iotas in homology}
 {f_{\iota_1}}_* = \left(\begin{array}{cc}
                         2 & 0 \\
                         0 & 1
                         \end{array}\right)\ ,\ 
                         &                     
 {f_{\iota_2}}_* = \left(\begin{array}{cc}
                         1 & 0 \\
                         0 & 2
                         \end{array}\right)\ ,\ 
                         &
 {f_{\iota_3}}_* = \left(\begin{array}{cc}
                         1 & -1 \\
                         1 & 1
                         \end{array}\right).
\end{eqnarray}

In Remark \ref{Remark: Meaning of the terms in the powers of z}, we said that intuitively
the terms depending on $z^i$ are associated to $i$-fold covers of the tori. In particular, the terms in $z^4$ are associated to double covers over
the $C_{\iota}$'s, which in turn can be identified with the elements of $H^1(C_{\iota},\Z/2)$.

The $4$-fold covers of $C$ that interest us are associated to pairs of double covers $(\iota_i, \iota_j)$ of $C$, with $\iota_i \neq \iota_j$.
Observe that the cover $\iota_i$ induces the map
$${f_{\iota_i}}_* : H_1(C_{\iota_i},\Z/2)\rightarrow H_1(C,\Z/2)$$
as defined in (\ref{Equation: The maps induced ba the iotas in homology}) (but with $\Z/2$ coefficients). Then the homomorphism
\begin{equation*}
 \iota_{j\circ i} := f_{\iota_i}^*(\iota_j) = \iota_j \circ {f_{\iota_i}}_* : H_1(C_{\iota_i},\Z/2) \rightarrow \Z/2
\end{equation*}
defines a class in $H^1(C_{\iota_i},\Z/2)$, and so a double cover $C_{\iota_{j\circ i}}$ of $C_{\iota_i}$ (cf. the map given in (5.22) of \cite{Ta1}).

\begin{Lemma} \label{Lemma: There is only one relevant 4-fold cover}
 For every pair of double covers $(\iota_i, \iota_j)$ of $C$ with $\iota_i \neq \iota_j$, the $4$-fold cover
 $$f_{\iota_i} \circ f_{\iota_{j\circ i}} : C_{\iota_{j\circ i}} \longrightarrow C$$
 is the unique determined by the image of the square $[0,1]^2 \subset \R^2$ by the matrix
 $$(l_{\iota_i} \circ l_{\iota_{j\circ i}})_* = \begin{pmatrix}
                                                 2 & 0 \\
                                                 0 & 2
                                                \end{pmatrix}$$
\end{Lemma}
\proof
We show only the case $(\iota_1,\iota_3)$, leaving to the reader the analogue computations for the other cases. First we have
\begin{equation*}
 \iota_{3\circ 1} = \iota_3 \circ {f_{\iota_1}}_* = \begin{pmatrix}
                                                 1 &  1
                                                \end{pmatrix}
                                                \begin{pmatrix}
                                                 2 & 0 \\
                                                 0 & 1
                                                \end{pmatrix} \sim 
                                                \begin{pmatrix}
                                                 0 &  1
                                                \end{pmatrix} \in H^1(C_{\iota_1},\Z/2).
\end{equation*}
Then $\iota_{3\circ 1}$ is the double cover $\iota_2$ of $C_{\iota_1}$ and:
\begin{equation*}
 (l_{\iota_1} \circ l_{\iota_{3\circ 1}})_* = {f_{\iota_1}}_* \circ {f_{\iota_2}}_* = 
                                                \begin{pmatrix}
                                                 2 & 0 \\
                                                 0 & 1
                                                \end{pmatrix}
                                                \begin{pmatrix}
                                                 1 & 0 \\
                                                 0 & 2
                                                \end{pmatrix} =
                                                \begin{pmatrix}
                                                 2 &  0 \\
                                                 0 & 2
                                                \end{pmatrix}.
\end{equation*}
\endproof
\noindent Observe in particular that the lemma implies that the $4$-fold covers of $C$ determined by $(\iota_i, \iota_j)$ and
$(\iota_j, \iota_i)$ are the same. We will denote by $\iota_4$ the $4$-fold cover of $C$ given by last lemma and by
$f_{\iota_4}:C_{\iota_4} \rightarrow C$ the corresponding covering map. 

When we defined $\GT_{W}$ we got a polynomial with formal variables $z_{B_W}$, where the $B_W$'s are homology classes in $W$ induced by homology
classes $B$ on the pull-back bundles over the holomorphic curves. In the case of tori it will be convenient to keep track of the projection on $X$ of
the classes $B_W$. 

\begin{Def} \label{Definition: GT_X for embedded tori}
 Let $C \subset X$ be an embedded torus. We define
% \begin{equation*}
%  \GT_X(\pi^*C) := \sum_{A \in H_2(X)} \mathrm{GT}_X(\pi^*C,A) t_{A}\ \in\ \Z[[H_2(W) \oplus H_2(X)]]
% \end{equation*} 
% where
% \begin{equation*}
%  \mathrm{GT}_X(\pi^*C,A) = \sum_{\substack{B \in H_2(\pi^*C)  \mbox{ \scriptsize{s.t.}}\\ \pi_*(B_{W})=A}}\mathrm{GT}(\pi^*C,\omega_C,B) z_{B_W}\ \in\ \Z[[H_2(W)]],
% \end{equation*}
 \begin{equation*}
  \mathrm{GT}_X(\pi^*C) := \sum_{\substack{B \in H_2(\pi^*C)}}\mathrm{GT}(\pi^*C,\omega_C,B) z_{B_W}t_{\pi_*(B_{W})}\ \in\ \Z[[H_2(W) \oplus H_2(X)]],
 \end{equation*}
 where $\omega_C$ is any symplectic form like in (\ref{Equation: Symplectic form above the curves}).
\end{Def}
 
Observe that the informations contained in $\mathrm{GT}_X(\pi^*C)$ are exactly the same as $\mathrm{GT}_W(\pi^*C)$ (as defined in last subsection),
which can be recovered just by setting $t_A = 1$ for all $A \in H_2(X)$. It is convenient to regard $\mathrm{GT}_X(\pi^*C)$
as a series in the variables $t_A$ and coefficients that are series in the variables $z_{B}$. Note finally that, by definition, if the coefficient
of $t_A$ in $\mathrm{GT}_X(\pi^*C)$ is non-zero, then $A$ is a non-negative multiple of $[C]$.

Now, given an embedded torus $C \subset X$ with $[C] \in \T$ and one of the four relevant covers $\iota :C_{\iota} \rightarrow C$ above, 
let us denote by $\pi^*C_{\iota}:= f_{\iota}^*(\pi^*C)$ the total space
of the bundle over $C_{\iota}$ and fiber $F$ obtained by pulling back via $f_{\iota}$ the bundle structure of $\pi^*C$. Then
the smooth $4$-manifold $\pi^*C_{\iota}$ is the total space of a ($2$- or $4$-fold) cover of $\pi^*C$ and we will call
$\pi^*f_{\iota}$ the corresponding covering map.
\begin{center} 
% \begin{tikzcd}[scale=2.5]
%  \pi^*C_{\iota} \arrow{r}{\pi^*f_{\iota}} \arrow{d}{\pi} & C_{\iota} \arrow{d}{\pi}\\
%  C_{\iota} \arrow{r}{f_{\iota}} & C
% \end{tikzcd}

\begin{tikzpicture}
  \matrix (m) [matrix of math nodes,row sep=3em,column sep=3em,minimum width=2em] {
     \pi^*C_{\iota} & \pi^*C \\
     C_{\iota} & C \\};
  \path[-stealth]
    (m-1-1) edge node [left] {$\pi$} (m-2-1)
            edge node [above] {$\pi^*f_{\iota}$} (m-1-2)
    (m-2-1.east|-m-2-2) edge node [above] {$f_{\iota}$} node [above] {} (m-2-2)
    (m-1-2) edge node [right] {$\pi$} (m-2-2);
  \label{Figure: ciccio}
%   \caption{Bifurcation.}
\end{tikzpicture}
\end{center}

Observe that, if $\omega_C$ is a symplectic form over $\pi^*C$ given by Theorem \ref{Theorem: Thurston}, the pull-back
$\omega_{C_\iota}:=(\pi^*f_{\iota})^*(\omega_C)$ is symplectic over $\pi^*C_{\iota}$ and makes the fibers symplectic.
 
\begin{Def} \label{Definition: GT_X for covers of embedded tori}
 Let $C \subset X$ be an embedded torus with $[C] \in \T$ and $\iota :C_{\iota} \rightarrow C$ one of the covers $\{\iota_1,\iota_2,\iota_3,\iota_4\}$. Define
% \begin{equation*}
%  \GT_X(\pi^*C_{\iota}) := \sum_{A \in H_2(X)} \mathrm{GT}_X(\pi^*C_{\iota},A) t_{A}\ \in\ \Z[[H_2(W) \oplus H_2(X)]]
% \end{equation*} 
% where
% \begin{equation*}
%  \GT_X(\pi^*C_{\iota},A) := \sum_{\substack{B \in H_2(\pi^*C_{\iota}) \mbox{ \scriptsize{s.t.}}\\ ((\pi^*f_{\iota})_*(B))_W = B_W, \\ \pi_*(B_W) = A} } \mathrm{GT}(\pi^*C_{\iota},\omega_{C_\iota},B) z_{B_W}\ \in\ \Z[[H_2(W)]]. 
% \end{equation*} 
 \begin{equation*}
  \mathrm{GT}_X(\pi^*C_{\iota}) := \sum_{\substack{B \in H_2(\pi^*C_{\iota})}}\mathrm{GT}(\pi^*C_{\iota},\omega_{C_\iota},B) z_{B_W}t_{\pi_*(B_{W})}\ \in\ \Z[[H_2(W) \oplus H_2(X)]],
 \end{equation*}
 where, with slight abuse of notation, we set $B_W=\big((\pi^*f_{\iota})_*(B)\big)_W$
 %where $\mathrm{GT}_W(\pi^*C_{\iota},B)$ is computed with respect to any symplectic form on $\pi^*C_{\iota}$ and normalized (like $\mathrm{GT}_W(\pi^*C,B)$ in Definition
 %\ref{Definition: GT_W for embedded surfaces}) using the $\mathrm{Spin}^c$-structure $((\pi^*f_{\iota})^*\circ \pi^*)(\mathfrak{s}_{\omega})$ over $\pi^*C_{\iota}$.
 %and, for any $B$, $z_{B_W}$ is a formal
 %variable encoding the homology class $B_W = ((\pi^*f_{\iota})_*(B))_W \in H_2(W)$ (cf. Notation \ref{Notation: Homology classes B_W}).
\end{Def}

In what follows, we will need to keep track not only of the type $(\epsilon(C),s) \in \{\pm 1\} \times \{0,1,2,3\}$ of a torus $C$ defined in the
end of Subsection \ref{Subsection: Definition of GT}, but also
of the particular $s$ double covers of $C$ having negative sign. We will say then that \textit{$C$ is of type $(\epsilon(C),I)$} if
$I \subseteq \{\iota_1,\iota_2,\iota_3\}$ is the set of the $s$ covers of $C$ that have negative sign.

\begin{Def} Given $(W,X,\pi,F)$, let $C$ be an embedded holomorphic torus of type $(\epsilon(C),I)$ with homology class in $\T$. We define the following 
Laurent series in the formal variables $t_A$, $A\in H_2(X)$ and coefficients that are series in the variables $z_B$, $B\in H_2(W)$:
\begin{equation*}
 P_{\pi}(C)  := \begin{cases}
                \left(\GT_X(\pi^*C)\right)^{\epsilon(C)} & \mbox{if } {I} = \emptyset\\
                              &   \\
                \left(\dfrac{\GT_X(\pi^*C)}{\GT_X(\pi^*C_{\iota'})}\right)^{\epsilon(C)} & \mbox{if } {I} = {\{\iota'\}} \\
                              &   \\
                \left(\dfrac{\GT_X(\pi^*C) GT_X(\pi^*(C_{\iota_4}))}{\GT_X(\pi^*C_{\iota'})\GT_X(\pi^*C_{\iota''})}\right)^{\epsilon(C)} & \mbox{if } {I} = {\{\iota',\iota''\}}\\
                              &   \\
                \left(\dfrac{\GT_X(\pi^*C) GT_X(\pi^*(C_{\iota_4}))}{\GT_X(\pi^*C_{\iota_1})\GT_X(\pi^*C_{\iota_2})\GT_X(\pi^*C_{\iota_3})}\right)^{\epsilon(C)} & \mbox{if } {I} = \{\iota_1,\iota_2,\iota_3\}\\
               \end{cases}
\end{equation*}
%\noindent where the sign of the exponent is the one given by the type of $C$. 
\end{Def}

\vspace{0.3 cm}

\begin{Def}
 Given $A \in \T$ primitive, we define
 \begin{equation} \label{Equation: Qu(B,m) twisted}
  \Qu_{\pi}(A,m,\omega,J) := \sum_{\{(C_k,m_k)\}}\prod_k r_{\pi}(C_k,m_k)
 \end{equation}
 \noindent where the sum is over the sets of couples $(C_k,m_k) \in \mathcal{M}_{c_k A}(X) \times \N^*$ with $c_k\geq 1$ and $\sum_k c_k m_k = m$
 and $r_{\pi}(C_k,m_k)$ is defined to be the coefficient of $t_{m_k[C_k]}$ in the formal power series $P_{\pi}(C_k)$.
\end{Def}
\noindent When $\omega$ and $J$ are understood we will write $\Qu_{\pi}(A,m)$ instead of $\Qu_{\pi}(A,m,\omega,J)$.

\begin{Rmk}
 Note that
 $\Qu_{\pi}(A,m)$ is the coefficient of $t_{mA}$ in
 \begin{equation*}
  \prod_{C \mbox{\scriptsize{ embedded }}|\ [C] \in \T} P_{\pi}(C).
 \end{equation*}

\end{Rmk}

\begin{comment}
\noindent where, for any cover $\iota$ of $C$, $\pi^*C_{\iota}$ is (the total space of) the pull-back of the bundle $\pi^*C$ over $C$ via the map
$f_{\iota}:C_{\iota} \rightarrow C$ and
\begin{equation*}
 \GT_W(\pi^*C_{\iota}) = \sum_{B \in H_2(\pi^*C_{\iota})} \mathrm{GT}(\pi^*C_{\iota},B) z_{B_W}  
\end{equation*}
where ${B_W}$ is the image of $B$ under the homomorphism induced in homology by $i \circ {f_{\iota}}$, with $i:C \rightarrow X$
\end{comment}

%Since we are now dealing with Laurent series, we are allowing also negative powers for the formal variables $z_{A}$, and we will use the identification
%$z_A^{-1} \sim {z}_{-A}$.

\subsection{Twisted $\mathbf{\GT}$} \label{Subsection: Twisted GT}  $ $ \\
\indent We can now define the $\pi$-twisted Gromov-Taubes invariants exactly like in the standard case.
For any $d \in \N$, fix $d/2$ generic points in $X$. 

\begin{Def} \label{Definition: Twisted GT}
 Given $(W,X,\pi,F)$, we define the $\pi$-twisted Gromov-Taubes invariant of $(X,\omega,J)$ by
 \begin{equation*}
  \mathrm{GT}_{\pi}(X,\omega,J) := \sum_{A \in H_2(X)} \mathrm{GT}_{\pi}(X,A)t_A\ \in\ \Z[[H_2(W) \oplus H_2(X)]]
 \end{equation*}
 where
 \begin{equation*} 
  \mathrm{GT}_{\pi}(X,A) = \sum_{(A_i,n_i) \in \mathcal{D}(A)} d_A! \left(\prod_{A_i \notin \T} \frac{\mathfrak{Ru}_{\pi}(A_i)^{n_i}}{d_{A_i}!n_i!} \cdot
  \prod_{A_i \in \T} \mathfrak{Qu}_{\pi}(A_i,n_i)\right)
 \end{equation*}
 and $\mathcal{D}(A)$ is the set of decompositions given in Definition \ref{Definition: Decompositions D of A}.
\end{Def}

At this point some words about the variables encoding the homology classes in $W$ and $X$ is due. Think for example to 
a holomorphic curve $L \subset \pi^*C$ with $g_{[L]} > 1$ for some embedded holomorphic curve $C \subset X$ with $g_{[C]} > 1$. Its total
contribution to
$\mathrm{GT}_{\pi}(X,\omega,J)$ is
\begin{equation}
 \epsilon(C)\epsilon(L)z_{[L]_W}t_{[C]},
\end{equation}
where $z_{[L]_W}$ and $t_{[C]}$ are in general independent (and in particular we could have $\pi_*([L]_W) \neq [C]$). On the other hand,
if $g_{[C]} = 1$, the total
contribution of $L$ is more complicated and in general it is not even of the form $ct_{[C]}$ for some Laurent series $c$ with variables $z_B$,
$B\in H_2(W)$. Roughly speaking the problem is that $L$ could ``appear'' also in some cover $\widetilde{C}$ of $\pi^*C$, in which
case we do not know on which $t_{[\widetilde{C}]}$ the contribution of $L$ should depend: the reason for which to set $t_{\pi_*{([L]_W)}}$ for the
$H_2(X)$-component of the weight of $L$ is a natural choice will appear more clear in Section \ref{Section: Mapping tori and twisted Lesfchetz zeta functions}.

\subsection{Proof of the symplectic invariance} \label{Subsection: Proof of the invariance}  $ $ \\
\indent In Subsection \ref{Subsection: About the symplectic invariance of GT} we gave a short overview towards the proof of the symplectic invariance
of the standard $\GT(X)$. Since the signs of the weights associated to the holomorphic curves are the same of Taubes, to prove the 
symplectic invariance of twisted $\GT$ it will be enough to carefully analyze the purely topological part of Taubes' proof. By this we mean the
following: given $(W,X,\pi,F)$ and a holomorphic curve $C \in \mathcal{M}_A(X,J,\Omega_{d_A})$, we want to
understand what happens to $\pi^*C$ and, eventually, its relevant covers, when we change $\Omega_{d_A}$, $J$ or we follow a smooth path of symplectic structures.
The answer is mostly completely contained in the proofs in \cite[Section 5]{Ta1} and \cite[Sections 4 and 5]{Ru}.
In what follows we will use the notations introduced in Subsection \ref{Subsection: About the symplectic invariance of GT}.

\begin{Lemma} \label{Lemma: M does not change by locally change J and Omega}
 For any $A \in H_2(X)$ and $(J,\Omega_d) \in \mathcal{U}_{d}$ (where $d=d_A$), there is an open neighborhood $U= U(J,\Omega_d)$
 of $(J,\Omega_d)$ in $\mathcal{U}_d$ such that for any $(J',\Omega_d') \in U$, there is an orientation preserving homeomorphism of
 ($0$-dimensional compact oriented) manifolds
 $$\mathcal{M}_A(X,J,\Omega_d) \cong \mathcal{M}_A(X,J',\Omega_d')$$
 that sends any $C \in \mathcal{M}_A(X,J,\Omega_d)$ to a $C' \in \mathcal{M}_A(X,J',\Omega_d')$ which is isotopic and close to $C$ in $X$.
\end{Lemma}

Lemma \ref{Lemma: M does not change by locally change J and Omega} is essentially the Assertion 2 of \cite[Proposition 5.2]{Ta1}.
The ``isotopic'' part is not explicitly stated but follows from Taubes' proof. Briefly, $\mathcal{M}_A(X,J,\Omega_d)$ 
is the counter-image of a regular value $(J,\Omega_d)$ of a certain smooth function with image in $\mathcal{A}_{d}$. By the Sard-Smale theorem,
the set of the regular values of this function must be an open and dense subset of $\mathcal{A}_{d}$ and the inverse function theorem implies
that if $(J',\Omega'_d)$ is close enough to $(J,\Omega_d)$ in $\mathcal{A}_{d}$, then there is a trivial cobordism from $\mathcal{M}_A(X,J,\Omega_d)$
to $\mathcal{M}_A(X,J',\Omega_d')$ that gives a homotopy in $X$ between each $C$ to a $C'$. For more details see the first steps of the proof of
\cite[Proposition 5.2]{Ta1} (and also \cite[Section 5]{Ru}).

The fact that we can actually get an isotopy comes by further restricting the choice for the couples $(J,\Omega_d)$.
Taubes observes in fact that, if $A$ is not multiply toroidal, the values of the time of the homotopy above where this fails to give an embedding
correspond to the degenerations appearing in the Gromov compactness theorem (\cite{Gr}). On the other hand, Taubes proves
(see step 5 of the proof of the aforementioned Taubes' proposition) that there exists an open and dense subset (i.e. $\mathcal{U}_{d}$) of $\mathcal{A}_{d}$
for which no such degeneration can occur, and the isotopy follows for $g_{[A]} > 1$. Finally, some local degenerations (like bifurcations)
that could a priori still happen if $A$ is multiply toroidal are excluded (for $U$ small enough) in the final steps of the proof of Taubes'
proposition.

\begin{Cor} \label{Corollary: Ru and Qu twisted are locally invariant}
 Let $(W,X,\pi,F)$ be a surface bundle, $A \in H_2(X)$ and $m\in \N$. Then $\Ru_{\pi}(A,\omega,J,\Omega_{d_A})$ and $\Qu_{\pi}(A,m,\omega,J)$
 do not depend on the choice of $(J,\Omega_{d_A})$ in a small neighborhood $U$ of $(J,\Omega_{d_A})$ in $\mathcal{A}_{d}$.
\end{Cor}
\proof
 If $C$ and $C'$ are embedded isotopic surfaces in $X$ with $g_{[C]} >1$ then the two $4$-manifolds $\pi^*C$ and $\pi^*C'$ are isotopic in $W$ and
 $\GT_W(\pi^*C) = \GT_W(\pi^*C')$. Observe that here we did not use only the fact that $\pi^*C$ and $\pi^*C'$ are diffeomorphic, but also the fact
 that (by Lemma \ref{Lemma: The acs of the Thurston symplectic forma are the same} and the statement of ``$\GT = \mathrm{SW}$''), $\GT_W(\pi^*C)$
 and $\GT_W(\pi^*C')$ have the same normalization of the variables.
 
 Similarly, if $g_{[C]} =1$, it is easy to check that $\GT_X(\pi^*C) = \GT_X(\pi^*C')$ and, identifying the $2$- and $4$-fold covers using the isotopy,
 for any $\iota$ we have also
 $\GT_X(\pi^*C_{\iota}) = \GT_X(\pi^*C'_{\iota})$. The result follows then from Lemma \ref{Lemma: M does not change by locally change J and Omega}.
\endproof

\begin{Cor} \label{Corollary: Ru twisted is a symplectic invariant}
 Let $(W,X,\pi,F)$ be a surface bundle, $A \in H_2(X) \setminus \T$  and $\{\omega_t\ |\ t\in [0,1]\}$ a smooth path of symplectic structures
 on $X$. Then 
 $$\Ru_{\pi}(A,\omega_0) = \Ru_{\pi}(A,\omega_1).$$
\end{Cor}
\proof
 The proof is very similar to that of last corollary and uses the smooth cobordism $\mathcal{C}_A(X,\{J_t\},\Omega_d)$ given by
 Lemma \ref{Lemma: A path of symplectic structures gives a smooth fibered manifold}. The main difference is that here the cobordism may be non-trivial
 since the projection $\mathcal{C}_A(X,\{J_t\},\Omega_d) \rightarrow [0,1]$ can have critical points (other degenerations are excluded by Taubes in Step
 7 of the proof of Proposition 5.2 of \cite{Ta1}). These critical points correspond
 to births or deaths of couples of holomorphic curves $C_1$ and $C_2$ with $\epsilon(C_1) = - \epsilon(C_2)$ whose total contribution, in the
 non-twisted case, obviously does not affect $\Ru$. Furthermore $C_1$ ad $C_2$ must be in the same connected component of
 $\mathcal{C}_A(X,\{J_t\},\Omega_d)$, and this gives an isotopy in $X$ between $C_1$ and $C_2$ so that their total contribution
 $$\epsilon(C_1)\GT_W(\pi^*C_1) +\epsilon(C_2)\GT_W(\pi^*C_2)$$
 to $\Ru_{\pi}$ is also zero, since $\GT_W(\pi^*C_1) = \GT_W(\pi^*C_2)$.
\endproof

\vspace{0.3 cm}

The proof of the symplectic invariance for $\Qu_{\pi}$ is more delicate since, for a given $B \in \T$, the $1$-dimensional manifold 
$\mathcal{C}_B(X,\{J_t\})$ is not in general compact. As explained in Subsection \ref{Subsection: About the symplectic invariance of GT},
a sequence of tori in $\mathcal{C}_{2B}(X,\{J_t\})$ could converge, for $t \rightarrow t_0$, to a $2$-fold cover of
a holomorphic torus in $\mathcal{M}_{B}(X,J_{t_0})$ (called the weak limit of the sequence).

Now, if $A\in \T$ is primitive and the path $\{\omega_t\ |\ t\in [0,1]\}$ is chosen carefully, then, for any $n \in \N$,
these weak limits correspond to the bifurcation points of $\mathcal{K}_{nA}(X,\{J_t\})$
(defined in (\ref{Equation: Union K of 1-dimensional spaces of tori given by a path of sympl. str.})) and these are exactly the points of
$\mathcal{K}_{nA}(X,\{J_t\})$ where the latter fails to be compact (recall that $\mathcal{K}_{nA}(X,\{J_t\})$ has the topology induced by the disjoint
union). Fix $n \in \N$ and consider the natural projection
\begin{equation*}
\begin{array}{ccccc}
 v_{nA} & : &\mathcal{K}_{nA}(X,\{J_t\}) & \longrightarrow & [0,1]\\
  & &(t,C)  & \longmapsto & t.
\end{array}
\end{equation*}

By \cite[Lemma 5.8]{Ta1}, $\mathcal{K}_{nA}(X,\{J_t\})$ is a $1$-dimensional manifold and we can assume that its bifurcation points and the
critical points
of $v_{nA}$ (births and deaths of pairs of tori with opposite signs) are finite and the corresponding values of $t \in [0,1]$ are all distinct and in 
$(0,1)$. Let us call $t \in (0,1)$ a \emph{bifurcation value} of $v_{nA}$ if $v_{nA}^{-1}(t_0)$ contains a bifurcation point and, with
a slight abuse of language, we will say that $t_0$ is a \emph{regular value} of $v_{nA}$ if it is neither a critical nor a bifurcation value.

By the aforementioned lemma of Taubes, if $t$ is a regular value of $v_{nA}$ and $C \in \mathcal{M}_{mA}(X,J_t)$ is an embedded torus for some
$m \leq n$, the signs of $C$ and its relevant covers are defined, and so we can compute the functions $P(C,\cdot)$ and $P_{\pi}(C)$. 

\begin{Prop}
 Let $(W,X,\pi,F)$ be a surface bundle and $\{\omega_t\ |\ t\in [0,1]\}$ a smooth path of symplectic structures on $X$. Let $\{J_t\ |\ t \in [0,1]\}$
 be a path of $\omega_t$-compatible almost complex structures given by Lemma
 \ref{Lemma: A path of symplectic structures gives a smooth fibered manifold}. Then for any $A \in \T$ primitive, if $t$ is a regular
 value of $v_{mA}$ for all $m$, the product
 \begin{equation*}
  \prod_{m > 0} \left( \prod_{C \in \mathcal{M}_{mA}(X,J_t)} P_{\pi}(C)\right)
 \end{equation*}

 \begin{comment}
  Then the product 
 \begin{equation} \label{Equation: Product of all the twisted P of C}
  \prod_{\substack{A \in \T \\ \mbox{\scriptsize{ primitive }}}} \left( \prod_{m \geq 1} \left( \prod_{C \in \mathcal{M}_{mA}(X,J_t)} P_{\pi}(C)\right)\right)
 \end{equation}
\end{comment}
 
\noindent does not depend on the choice of $t$.
\end{Prop}
\proof
 Let $t_0$ be either a critical or a bifurcation value of the map $v_{nA}$ for some $n \in \N$ and let $r >0$ such that $t_0$ is the only
 non-regular value of $v_{nA}$ in the interval $(t_0-r,t_0+r)$. We want to prove that for any $0<s<r$
 \begin{equation} \label{Equation: The two product of tori are the same before and after a critical point}
  \prod_{m \leq n} \left( \prod_{C \in \mathcal{M}_{mA}(X,J_{t_0-s})} P_{\pi}(C)\right) =  \prod_{m \leq n} \left( \prod_{C \in \mathcal{M}_{mA}(X,J_{t_0+s})} P_{\pi}(C)\right).
 \end{equation}
 Note first that, since $(t_0-r,t_0) \cup (t_0,t_0+r)$ contains only regular values of $v_{nA}$, the two sides of
 (\ref{Equation: The two product of tori are the same before and after a critical point}) do not depend on the choice of $s$.

 Let us first check the case where $t_0$ corresponds to a critical point. As in the proof of Corollary \ref{Corollary: Ru twisted is a symplectic invariant}, $t_0$
 corresponds to the birth or a death of a pair of isotopic embedded tori $C$ and $C'$ both in $v_{nA}^{-1}(t_0+s)$ or,
 respectively, $v_{nA}^{-1}(t_0-s)$. By \cite[Lemma 5.9]{Ta1} 
 \begin{eqnarray*}
  \epsilon(C) = - \epsilon(C') & \mbox{ and } & \epsilon(C_{\iota}) = \epsilon(C'_{\iota})
 \end{eqnarray*}
 where $\iota$ is any of the connected covers $\{\iota_1,\iota_2,\iota_3,\iota_4\}$, so that if $C$ is of type $(\epsilon,I)$ then $C'$
 is of type $(-\epsilon,I)$.
 Moreover, since $C$ and $C'$ are isotopic, for any $\iota$ we have 
 %(and identifying the (co)homology classes of $C$ and $C'$ using the isotopy),
 % $\pi^*(C_{\iota})$ and $\pi^*(C'_{\iota})$ are diffeomorphic and
 $\GT_X(\pi^*C_{\iota}) = \GT_X(\pi^*C'_{\iota}).$
 Summarizing: 
 $$P_{\pi}(C') = (P_{\pi}(C))^{-1}$$
 and the total contribution of $C$ and $C'$ to both sides of (\ref{Equation: The two product of tori are the same before and after a critical point})
 is
 $$P_{\pi}(C) P_{\pi}(C') = 1.$$
 
 Let us prove now the equivalence in (\ref{Equation: The two product of tori are the same before and after a critical point}) for $t_0$ corresponding
 to a bifurcation point. This means that there exists $B = kA$ for some $k$ and a sequence of tori in $\mathcal{C}_{2B}(X,\{J_t\})$ that converges to
 a double cover $\overline{\iota}$ of a torus $C_0 \in \mathcal{M}_{B}(X,J_{t_0}) \subset \mathcal{C}_{B}(X,\{J_t\})$. There are two possibilities,
 depending on whether the sequence converges to $C_0$ for increasing or decreasing $t$. We will assume $t$ decreasing and we will leave the (completely
 analogue) proof for the other case to the reader.

 \begin{figure} [ht] 
  \begin{center}
   \includegraphics[scale = .63]{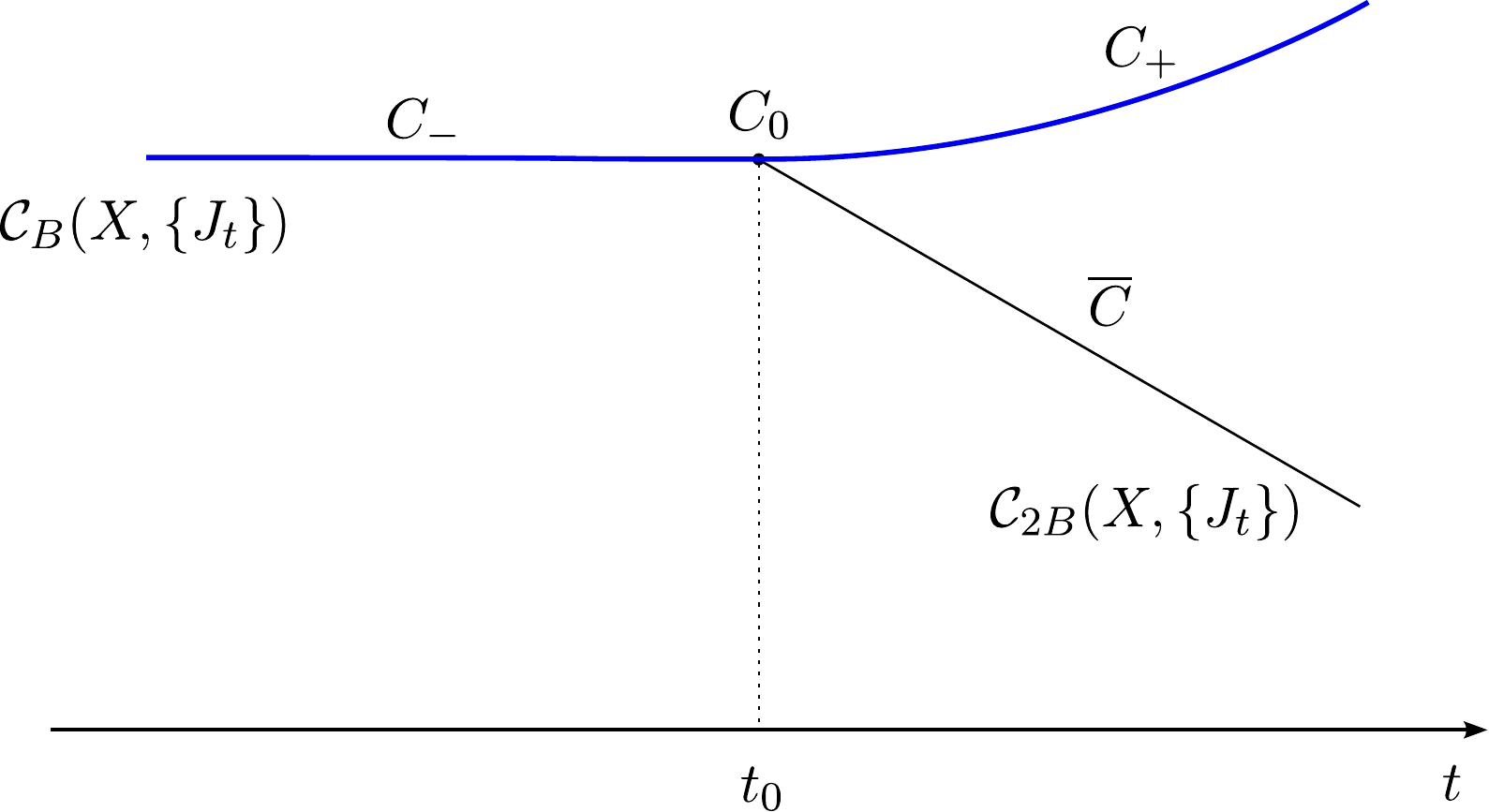}
   \caption{Bifurcation.}
   \label{Figure: Bifurcation}
  \end{center}
 \end{figure}

 The connected component of $\mathcal{C}_{B}(X,\{J_t\})$ containing $C_0$ gives an isotopy between the two embedded tori
 $C_- \subset \mathcal{M}_{B}(X,J_{t_0-s})$ and $C_+ \subset \mathcal{M}_{B}(X,J_{t_0+s})$, so that
 \begin{equation} \label{Equation: GT_X of C_- C_+ and C_0 equals}
  \GT_X(\pi^*C_-) = \GT_X(\pi^*C_0) = \GT_X(\pi^*C_+)
 \end{equation}
 and, for any $\iota \in \{\iota_1, \iota_2, \iota_3, \iota_4\}$,
 \begin{equation} \label{Equation: GT_X of covers of C_- C_+ and C_0 equals}
  \GT_X(\pi^*(C_-)_{\iota}) = \GT_X(\pi^*(C_0)_{\iota}) = \GT_X(\pi^*(C_+)_{\iota}),
 \end{equation}
 where, as usual, we identified the covers of different tori using the isotopy. By Lemma 5.10 of \cite{Ta1} we can recover the relative signs of
 $C_-$, $C_+$ and their double covers:
  \begin{equation} \label{Equation: Signs of C_- in terms of C_+ for the bifurcation}
  \begin{array}{lll}
   \epsilon(C_-) & = & \epsilon(C_+); \\
      & & \\
   \epsilon((C_-)_{\overline{\iota}}) & = & - \epsilon((C_+)_{\overline{\iota}});  \\
      & & \\
   \epsilon((C_-)_{\iota}) & = & \epsilon((C_+)_{\iota})  \mbox{ for any double cover } \iota \neq \overline{\iota}.
  \end{array}
 \end{equation}
 
 Now, let $\overline{C}$ be the only torus in $\mathcal{M}_{2B}(X,J_{t_0+s})$ belonging to the relevant connected component of
 $\mathcal{C}_{2B}(X,\{J_t\})$. This connected component gives a homotopy between $\overline{C}$ and the double cover $(C_0)_{\overline{\iota}}$
 (and so also with the double covers $(C_-)_{\overline{\iota}}$ and $(C_+)_{\overline{\iota}}$). Then, by homotopy invariance and
 (\ref{Equation: GT_X of covers of C_- C_+ and C_0 equals}), we have
  \begin{equation} \label{Equation: GT_X of C-bar and iota-bar covers of C_- and C_+ equals}
  \GT_X(\pi^*\overline{C}) = \GT_X(\pi^*((C_-)_{\overline{\iota}})) = \GT_X(\pi^*((C_+)_{\overline{\iota}})).
 \end{equation}
 
% Let $\iota \neq \overline{\iota}$ be any of the two remaining non-trivial double covers of $C_0$ and consider the non trivial double cover
% $\overline{C}_{\widehat{\iota}}$ where
% \begin{equation} \label{Equation: iota widehat}
%  \widehat{\iota} := f_{\overline{\iota}}^*(\iota).
% \end{equation}
 
 Moreover, observe that the $4$-fold cover $\iota_4$ of $C_0$ given by Lemma \ref{Lemma: There is only one relevant 4-fold cover}
 naturally induces a double cover $\widehat{\iota}$ of $\overline{C}$. More in detail, as for $s \rightarrow 0\ $ $\overline{C}$ converges to (the 
 total space of) the double cover $\overline{\iota}$ of $C_0$, we require that, for $s \rightarrow 0$, $\overline{C}_{\widehat{\iota}}$ converges to
 (the total space of) the double cover $f_{\overline{\iota}}^*(\iota)$ of $(C_0)_{\overline{\iota}}$ for any (non-trivial) $\iota \neq \overline{\iota}$.
 Since $f_{\overline{\iota}}^*(\iota)$ is exactly the cover $\iota_4$ of $C_0$, again by isotopy invariance, we have:
 \begin{equation}
  \GT_X(\pi^*(\overline{C}_{\widehat{\iota}})) = \GT_X(\pi^*((C_-)_{\iota_4})) = \GT_X(\pi^*((C_+)_{\iota_4})).
 \end{equation}

% be the unique double cover of $\overline{C}$ such that the map
% $$f_{\widehat{\iota}}: T^2 = \overline{C}_{\widehat{\iota}} \longrightarrow \overline{C} \subset X$$
% is homotopic to
% $$f_{\iota_4}: T^2 = \overline{C}_{\iota_4} \longrightarrow C_0 \subset X.$$

 By Lemma 5.11 of \cite{Ta1} we can recover the signs of $\overline{C}$ and of its double covers:
 \begin{equation} \label{Equation: Signs of overline C in terms of C_+ for the bifurcation}
  \begin{array}{lll}
   \epsilon(\overline{C}) & = & -\epsilon(C_+) \epsilon((C_+)_{\overline{\iota}});\\
    & & \\
   \epsilon(\overline{C}_{\widehat{\iota}}) & = & \prod_{\iota \neq \overline{\iota}}\ \epsilon((C_+)_{{\iota}});\\
      & & \\
   \epsilon(\overline{C}_{\iota}) & = & +1 \mbox{ for } \iota \neq \widehat{\iota}.
  \end{array}
 \end{equation}

 We are now ready to prove that the equivalence (\ref{Equation: The two product of tori are the same before and after a critical point}) 
 holds also when we cross a bifurcation. To prove the result it is enough to show that 
 $$ P_{\pi}(C_+)\cdot P_{\pi}(\overline{C}) = P_{\pi}(C_-).$$ 
 We will check only the three cases below, corresponding to the three possible types $(+,I)$ of $C_-$ with $\epsilon(C_-) = +$
 and $\overline{\iota} \notin I$: if $\epsilon(C_-) = -$ each case will be the reciprocal of the corresponding case for $\epsilon(C_-) = +$,
 while the cases with $\overline{\iota} \in I$ do not give new relations and are left to the reader. 
 %In what follows $\overline{\iota}$
 %will still indicate the double cover of $C_0$ which appears as weak limit of $\overline{C}$ for $t \rightarrow t_0$.
 
 \vspace{0.3 cm}
 
 \noindent \textbullet $\ [I = \emptyset]$: $C_+$ is of type $(+,\{\overline{\iota}\})$ and $\overline{C}$ is of type $(+,\emptyset)$. Then:
 \begin{eqnarray*}
 P_{\pi}(C_+)\cdot P_{\pi}(\overline{C}) & = & \dfrac{\GT_X(\pi^*C_+)}{\GT_X(\pi^*((C_+)_{\overline{\iota}})}\cdot \GT_X(\pi^*\overline{C}) \\
                                         &   & \\
                                         & = & \GT_X(\pi^*C_+)  \\
                                         &   & \\
                                         & = & \GT_X(\pi^*C_-) = P_{\pi}(C_-).
 \end{eqnarray*}
  
  \vspace{0.3 cm}
  
 \noindent \textbullet $\ [I=\{\iota'\}]$: $C_+$ is of type $(+,\{\overline{\iota},\iota'\})$ and $\overline{C}$ is of type $(+,\{\widehat{\iota}\})$,
 where $\widehat{\iota}$ is as above.\\ \noindent \phantom{c} Then:
  
 \begin{eqnarray*}
  P_{\pi}(C_+)\cdot P_{\pi}(\overline{C}) & = & \dfrac{\GT_X(\pi^*C_+) GT_X(\pi^*(({C_+})_{\iota_4}))}{\GT_X(\pi^*((C_+)_{\overline{\iota}}))\GT_X(\pi^*((C_+)_{\iota'}))} \cdot \dfrac{\GT_X(\pi^*\overline{C})}{\GT_X(\pi^*\overline{C}_{\widehat{\iota}})}\\
                                          &   & \\
                                          & = & \dfrac{\GT_X(\pi^*C_+)}{\GT_X(\pi^*((C_+)_{\iota'}))}\\
                                          &   & \\
                                          & = & \dfrac{\GT_X(\pi^*C_-)}{\GT_X(\pi^*((C_-)_{\iota'}))} =  P_{\pi}(C_-).                                      
 \end{eqnarray*}

  \vspace{0.3 cm}  
     
 \noindent \textbullet $\ [I=\{\iota',\iota''\}]$: $C_+$ is of type $(+,\{\overline{\iota},\iota',\iota''\})$ and $\overline{C}$ is of type
  $(+,\emptyset)$ (since by (\ref{Equation: Signs of overline C in terms of C_+ for the bifurcation})\\ \noindent \phantom{c} we have 
  $\epsilon(\overline{C}_{\widehat{\iota}}) = \epsilon((C_+)_{{\iota'}})\cdot \epsilon((C_+)_{{\iota''}}) = (-1) \cdot (-1) = +1$). Then:
  \begin{eqnarray*}
   P_{\pi}(C_+)\cdot P_{\pi}(\overline{C}) & = & \dfrac{\GT_X(\pi^*C_+) GT_X(\pi^*(({C_+})_{\iota_4})) \cdot \GT_X(\pi^*\overline{C})}{\GT_X(\pi^*((C_+)_{\overline{\iota}}))\GT_X(\pi^*((C_+)_{\iota'}))\GT_X(\pi^*((C_+)_{\iota''}))} \\
   %                                        &   & \\
   %                                        & = & \dfrac{\GT_X(\pi^*C_+) GT_X(\pi^*(({C_+})_{\iota_4})) \cdot \GT_X(\pi^*((C_+)_{\overline{\iota}}))}{\GT_X(\pi^*((C_+)_{\overline{\iota}}))\GT_X(\pi^*((C_+)_{\iota'}))\GT_X(\pi^*((C_+)_{\iota''}))} \\
                                           &   & \\
                                           & = & \dfrac{\GT_X(\pi^*C_+) GT_X(\pi^*(({C_+})_{\iota_4}))}{\GT_X(\pi^*((C_+)_{\iota'}))\GT_X(\pi^*((C_+)_{\iota''}))} \\
                                           &   & \\
                                           & = & \dfrac{\GT_X(\pi^*C_-) GT_X(\pi^*(({C_-})_{\iota_4}))}{\GT_X(\pi^*((C_-)_{\iota'}))\GT_X(\pi^*((C_-)_{\iota''}))} =  P_{\pi}(C_-).
 \end{eqnarray*}

 \endproof

\begin{Cor} \label{Corollary: Qu twisted is a symplectic invariant}
 Let $(W,X,\pi,F)$ be a surface bundle, $A \in \T$  and $\{\omega_t\ |\ t\in [0,1]\}$ a smooth path of symplectic structures
 on $X$. Then 
 $$\Qu_{\pi}(A,m,\omega_0) = \Qu_{\pi}(A,m,\omega_1).$$
\end{Cor}

Theorem \ref{Theorem: GT twisted is a symplectic invariant} follows then from the definition of $\GT_{\pi}(X,\omega,J)$ and
corollaries \ref{Corollary: Ru and Qu twisted are locally invariant}, \ref{Corollary: Ru twisted is a symplectic invariant} and
\ref{Corollary: Qu twisted is a symplectic invariant}.

\section{Mapping tori and twisted Lesfchetz zeta functions} \label{Section: Mapping tori and twisted Lesfchetz zeta functions}

In this section we first review the definition of standard and twisted Lefschetz zeta functions of surface diffeomorphisms. We then define
the ``bundled twistings'' for Lefschetz zeta functions and prove the equivalence with certain Jiang's twisted Lefschetz zeta functions.
Finally we compute the twisted Gromov-Taubes invariants for symplectic closed $4$-manifolds of the form $S^1 \times Y$ and give
a proof of Theorem \ref{Theorem: GT twisted = Zeta twisted in intro} of introduction.

\subsection{Commutative Lefschetz zeta functions}  \label{Subsection: Commutative Lefschetz zeta functions}$ $ \\
\indent Let $S$ be an oriented compact connected surface and let $\phi : S \rightarrow S$ be an orientation preserving diffeomorphism.
A point $x \in S$ is a \emph{fixed point} of $\phi$ if $\phi(x) = x$. More in general, given $n \in \N^*$, $x$ is a \emph{periodic point of period $n$}
(or \emph{$n$-periodic point}) of $\phi$ if $\phi^n(x) = x$. An $n$-periodic point $x$ is said to be \emph{non-degenerate} if 
$$\det(\mathbbm{1} - d_x\phi^n) \neq 0.$$
The diffeomorphism $\phi$ is \emph{non-degenerate} if for every $n \in \N^*$ there is only a finite number of $n$-periodic points, each of which is non-degenerate.
Note in particular that in this case an $n$-periodic point $x$ is, for every $p \in \N^*$, also a $pn$-periodic point, which is non-degenerate
with respect to $\phi^{pn}$. 
It is well known that every element of the oriented mapping class group $\mathrm{MCG}(S)$ has a representative that is smooth and non degenerate.
From now on we will then assume the non-degeneracy of $\phi$.

The \emph{Lefschetz sign} of an $n$-periodic point $x$ is 
\begin{equation*}
 \varepsilon(x,n) := \mathrm{sign}(\det(\mathbbm{1} - d_x\phi^n)).
\end{equation*}
Observe that any $n$-periodic point of $\phi$ can be naturally interpreted as a fixed point of $\phi^n$ and that the Lefschetz signs given by the
two interpretations are the same. On the other hand, we remark that in general $(\varepsilon(x,n))^p \neq \varepsilon(x,pn)$ for $p >1$.

Let now
$$\To_{\phi} := \dfrac{S \times [0,+\infty)}{(x,t+1)\sim (\phi(x),t)}$$
be the \emph{mapping torus of $(S,\phi)$}. 
%, if $S$ and $\phi$ are understood, we will write just $\To$ instead of $\To_{\phi}$. Moreover
We will often use the identification $S = S\times \{0\} \subset \To_{\phi}$.

Let $t$ be the coordinate of $[0,+ \infty)$ in $\To_{\phi}$. An \emph{$n$-periodic orbit of $\To_{\phi}$} is a closed positively oriented orbit $\delta$
of the flow on $\To_{\phi}$ of the vector field $\partial_t$, such that $\langle \delta,S \rangle = n$ and considered up to reparametrization.
An orbit is \emph{simple} if it is not a non-trivial cover of another closed orbit of the flow. Any orbit can then be identified with some $p$-fold
cover $\delta^p$ of some simple orbit $\delta$.

To any $n$-periodic simple orbit $\delta$ corresponds in a natural way the set $\{\delta \cap S\}$ of $n$ periodic points,
each of period $n$: the \emph{Lefschetz sign of $\delta$} is defined by
$$\varepsilon(\delta):= \varepsilon(x,n) \mbox{ for any } x\in \{\delta \cap S\}.$$
In general, given an integer $p>0$, the Lefschetz sign of $\delta^p$ is
$$\varepsilon(\delta^p):= \varepsilon(x,pn) \mbox{ for any } x\in \{\delta \cap S\}.$$

\begin{Def}
 Given a simple orbit $\delta$ of $(S,\phi)$, the \emph{local Lefschetz zeta function} of $\delta$ is the formal power series
 \begin{equation*}
  \zeta_{\delta}(t) := \exp\left(\sum_{m\geq 1} \varepsilon(\delta^m)\frac{t^m}{m} \right) \in \Z[[t]].
 \end{equation*}
\end{Def}

\begin{Rmk} 
 For the following, we remark that it is possible to give more explicit formulas for $\zeta_{\delta}(t)$ depending on the \emph{type} of $\delta$.
 It is a standard result that there are essentially two types of non-degenerate periodic orbits: \emph{elliptic} and \emph{hyperbolic}. A simple
 $n$-periodic orbit $\delta$ is elliptic if, for any $x \in \delta \cap S$, the eigenvalues of $d_x\phi^n$ are complex conjugates and is hyperbolic
 if they are real.  In the last case we can make a further distinction: $\delta$ is called \emph{positive} (resp. \emph{negative}) hyperbolic if the
 eigenvalues of $d_x\phi^n$ are both positive (resp. negative). It is then easy to prove that: 
 \begin{equation} \label{Equation: Lefschetz signs of the type of the orbits}
  \varepsilon(\delta^m) = \left\{ \begin{array}{ll}
                       +1 & \mbox{if } \delta \mbox{ elliptic;}\\
                       -1  & \mbox{if } \delta \mbox{ positive hyperbolic;}\\
                       (-1)^{m+1}  & \mbox{if } \delta \mbox{ negative hyperbolic}
                      \end{array}\right.
 \end{equation}
 and the following expressions follow (see for example Lemma 3.15 of \cite{Sp}):
 \begin{equation} \label{Equation: Local Lefschetz Zeta functions of the type of orbits}
  \zeta_{\delta}(t) = \left\{ \begin{array}{ll}
                       \frac{1}{1-t} = 1+t+t^2+\ldots & \mbox{if } \delta \mbox{ elliptic;}\\
                       1-t  & \mbox{if } \delta \mbox{ positive hyperbolic;}\\
                       1+t  & \mbox{if } \delta \mbox{ negative hyperbolic.}
                      \end{array}\right.
 \end{equation}
\end{Rmk}

\begin{Def} \label{Definition: Lefschetz zeta function standard}
 The \emph{Lefschetz zeta function} of $\phi$ is
 \begin{equation*}
  \zeta(\phi) := \prod_{\delta \mbox{\scriptsize{ simple}}} \zeta_{\delta}\left(t^{\langle \delta,S \rangle}\right) \in \Z[[t]]
 \end{equation*}
 where the product is meant to be over all the simple periodic orbits $\delta$ of $\phi$.
\end{Def}

We observe that in $\zeta(\phi)$ any periodic orbit of $\phi$ is counted exactly once and with weight depending only on its sign and its homology
class in $\To_{\phi}$.
The Lefschetz zeta function can be expressed also in terms of the \emph{Lefschetz numbers} of the iterates of $\phi$.
\begin{Def}
 Given $(S,\phi)$, for $i=0,1,2$ let ${\phi}_i: H_i(S) \rightarrow H_i(S)$ be the induced isomorphisms in homology.
 \emph{The Lefschetz number of $\phi$} is:
 $$\Le(\phi) := \sum_{i=0}^2 (-1)^i \Tr({\phi}_i).$$
\end{Def}

\begin{LFPThm}Let $\mathrm{Fix}(\phi)$ denote the set of fixed points of $\phi$. Then:
 $$\Le(\phi) = \displaystyle \sum_{x\in \mathrm{Fix}(\phi)} \varepsilon(x,1).$$
\end{LFPThm}
%\noindent The following is a classical result (see for example \cite{Fried}).
\begin{Cor} \label{Corollary: Lefschetz zeta function in terms of Lefschetz numbers}
 $ \zeta(\phi)= \exp\left( \displaystyle\sum_{n\geq 1} \Le(\phi^n)\dfrac{t^n}{n} \right).$
\end{Cor}
The proof of last corollary is matter of express any $n$-periodic orbit of $\phi$ as a fixed point
of $\phi^n$, reorganize the product in the definition of $\zeta(\phi)$ as a product in $n\geq 1$ and then apply the Lefschetz fixed
point theorem (see for example \cite{Fried}).

The last corollary implies in particular that $\zeta(\phi)$ is a topological invariant of $\To_{\phi}$. Moreover we have the following corollary (see
for example \cite{A' Campo}).
\begin{Cor} \label{Corollary: Zeta functions equals torsion in the Abelian case} Let $\tau\left(\To_{\phi}\right)=\tau\left(\To_{\phi},t\right)$ denote the Reidemeister torsion
 of $\To_{\phi}$. Then:
 $$\zeta(\phi) = \displaystyle\prod_{i=0}^2 \left(\det\left(\mathbbm{1} - t{\phi}_i\right)\right)^{(-1)^{i+1}} \doteq \tau\left(\To_{\phi}\right)$$
 where $\doteq$ denotes the equivalence up to multiplications by monomials of the form $\pm t^m$.
\end{Cor}
If not stated otherwise, all the Reidemeister torsions of mapping tori will be considered with the normalization given by the
co-oriented $0$-page $S$, so that the equivalence $\doteq$ with the corresponding zeta functions will be in fact an equality.

\vspace{0.3 cm}

There are various refinements of $\zeta(\phi)$: the richest Abelian that one can define is obtained by twisting the
contribution of each orbit by its total homology class in $\To_{\phi}$. By ``Abelian'' we mean that it can be obtained by using an Abelian
representation of $\pi_1(T_{\phi})$. This refinement is obtained by encoding the elements
$a \in H_1(\To_{\phi})$ in formal variables $t_a$ such that $t_a\cdot t_b = t_{a+b}$ like those in the definition of $\GT$ and
then considering the \emph{total Abelian Lefschetz zeta function}:
\begin{equation} \label{Equation: Total Abelian Lefschetz zeta function}
 \zeta_{\mathrm{A}}(\phi) := \prod_{\delta \mbox{\scriptsize{ simple}}} \zeta_{\delta}\left(t_{[\delta]} \cdot t^{\langle \delta,S \rangle}\right) \in \Z[[H_1\left(\To_{\phi}\right)]][[t]].
\end{equation}
The series $\zeta_{\mathrm{A}}(\phi)$ should be thought of as a series in the formal variable $t$ and coefficients in $\Z[[H_1\left(\To_{\phi}\right)]]$.

%Note that even if the information given by (the exponent of) the variable $t$ is contained by that given by $t_{[\delta]}$, it will be convenient to keep it in evidence.
All the results above about $\zeta(\phi)$ can be generalized to $\zeta_{\mathrm{A}}(\phi)$ by introducing corresponding
twisting coefficients for the Lefschetz numbers and the Reidemeister torsion (see \cite{Fried}, \cite{Fried 2} or \cite{Geo-Nic}). Briefly, instead of
considering the homomorphisms induced by $\phi$ in homology, one considers the homeomorphism 
$$\widetilde{\phi} : \widetilde{S} \longrightarrow \widetilde{S}$$
induced by $\phi$ on the universal Abelian cover $\widetilde{S}$ of $S$. Given a cellular decomposition of $\widetilde{S}$ obtained by lifting
a cellular decomposition of $S$, we can assume that both $\phi$ and $\widetilde{\phi}$ are cellular maps. Then, choosing basis and taking local
$\Z[H_1(S)]$-coefficients for $H_i(\widetilde{S})$, $i=0,1,2$, $\widetilde{\phi}$ gives rise to non-singular
$\Z[\mathrm{coker}(\mathbbm{1} - {\phi}_1)]$-matrices $\widetilde{\phi}_i$ of order $\mathrm{rk}(H_i(S))$.
The Mayer-Vietoris exact sequence implies that
\begin{equation} \label{Equation: Homology of a mapping torus}
 H_1\left(\To_{\phi}\right) \cong \mathrm{coker}(\mathbbm{1} - {\phi}_1) \oplus \Z_{\mu}
\end{equation}
so that the matrices ${\widetilde{\phi}}_i$ can be interpreted as matrices with coefficients in $\Z[H_1\left(\To_{\phi}\right)]$.
The \emph{total Abelian Lefschetz number of $\phi$} is 
\begin{equation} \label{Equation: Total Abelian Lefschetz numbers}
 \Le_{\mathrm{A}}(\phi) := \sum_{i=0}^2 (-1)^i \Tr(t_{\mu}\widetilde{\phi}_i) \in \Z[H_1\left(\To_{\phi}\right)].
\end{equation}
where $t_{\mu}$ is a formal variable encoding the homology class of the generator $\mu$ of $H_1\left(\To_{\phi}\right)$.
The \emph{total Abelian Lefschetz fixed point Theorem} (cf. Theorem 1 of \cite{Fried 2}) gives then
\begin{equation} \label{Equation: Total Abelian Lefschetz fixed point theorem}
 \Le_{\mathrm{A}}(\phi) = \sum_{x \in \mathrm{Fix}(\phi)}\varepsilon(x,1) [\delta_x] \in \Z[H_1\left(\To_{\phi}\right)]
\end{equation}
where $\delta_x$ is the unique $1$-periodic orbit of $\phi$ with $\delta \cap S = \{x\}$. We have then the following analogue of Corollary
\ref{Corollary: Lefschetz zeta function in terms of Lefschetz numbers}:
\begin{equation} \label{Equation: Total Abelian Lefschetz zeta function in terms of Lefschetz numbers}
 \zeta_{\mathrm{A}}(\phi)= \exp\left( \displaystyle\sum_{n\geq 1} \Le_{\mathrm{A}}(\phi^n)\dfrac{t^n}{n} \right) \in \Z[[H_1\left(\To_{\phi}\right)]].
\end{equation}
Similarly, the generalization to the total Abelian case of Corollary \ref{Corollary: Zeta functions equals torsion in the Abelian case} gives:
\begin{equation} \label{Equation: Total Abelian zeta functions and Reidemeister torsions}
 \zeta_{\mathrm{A}}(\phi) = \displaystyle\prod_{i=0}^2 \left(\det\left(\mathbbm{1} - t(t_{\mu}\widetilde{\phi}_i\right)\right)^{(-1)^{i+1}} = \tau_{\mathrm{A}}\left(\To_{\phi},t\right).
\end{equation}
where $\tau_{\mathrm{A}}\left(\To_{\phi},t\right)$ is the \emph{total Abelian Reidemeister torsion of $\To_{\phi}$}.

\begin{Ex} Let $L=K_1 \sqcup \ldots \sqcup K_n$ be a fibered $n$-component link in a homology $3$-sphere $Y$ and let $(S,\phi)$ be the corresponding
open book decomposition of $Y$ (so that in particular $\To_{\phi} \cong Y \setminus L$). Let $t_i,\ i\in\{1,\ldots,n\}$ be formal variables encoding the
homology classes of the meridians of the $K_i$'s. Then we can express all the variables $t_a,\ a\in H_1(\To_{\phi})$ in the definition of
 $\zeta_{\mathrm{A}}(\phi)$ in terms of the $t_i$'s and:
 \begin{equation}
 \zeta_{\mathrm{A}}(\phi)|_{t=1} = \tau_{\mathrm{A}}(Y\setminus L) =  \begin{cases}
                                                                                     \dfrac{\Delta_L(t_1)}{1-t_1} & \mbox{ if } n=1\\
                                                                                       &\\
                                                                                     \Delta_L(t_1,\ldots,t_n) & \mbox{ if } n>1
                                                                                    \end{cases}
 \end{equation}
 where $\zeta_{\mathrm{A}}(\phi)|_{t=1}$ is the evaluation in $t=1$ of the series $\zeta_{\mathrm{A}}(\phi)$ and $\Delta_L$ denotes the
 \emph{multivariable Alexander polynomial of $L$}. 
\end{Ex}

\subsection{Non-commutative Lefschetz zeta functions}  \label{Subsection: Non-commutative Lefschetz zeta functions}$ $ \\
\indent In this subsection we recall the \emph{twisted Lefschetz zeta functions} defined by Jiang in \cite{J} (see \cite{J-W}). Let $(S,\phi)$ be
as before, let $R$ be an Abelian ring with unity and $n \in \N^*$. Choose a base point $y_0 \in S \subset \To_{\phi}$ and consider a representation
\begin{equation*}
 \rho : \pi_1(\To_{\phi}) \longrightarrow \mathrm{GL}(n,R)
\end{equation*}
of the fundamental group $\pi_1(\To_{\phi})$ (for simplicity in the notation we will avoid to refer to $y_0$) into the linear group over $R$
of order $n$. Then $\rho$ extends to a representation (still denoted $\rho$)
\begin{equation*}
 \rho : \Z[\pi_1(\To_{\phi})] \longrightarrow \mathrm{M}(n,R)
\end{equation*}
of the group ring $\Z[\pi_1(\To_{\phi})]$ into the algebra of $n \times n$ matrices over $R$.

\vspace{0.3 cm}

For any simple orbit $\delta$ of $\To_{\phi}$, choose a path $\gamma_{\delta}$ from $y_0$ to any point
of $\delta$. Then, for any $m \geq 0$, $\gamma_{\delta} \ast \delta^m \ast \gamma_{\delta}^{-1}$
can be regarded as a closed path based on $y_0$.

\begin{Def}
 Given a simple orbit $\delta$ of $(S,\phi)$, the \emph{local $\rho$-twisted Lefschetz zeta function of $\delta$} is the formal power series
 \begin{equation*}
  \zeta_{\delta}(\rho,t) := \exp\left(\sum_{m\geq 1} \varepsilon(\delta^m) \Tr\left(\rho\left([\gamma_{\delta} \ast \delta^m \ast \gamma_{\delta}^{-1}]\right)\right) \frac{t^m}{m} \right) \in R[[t]].
 \end{equation*}
\end{Def}

Observe that $\zeta_{\delta}(\rho,t)$ does not depend on the choice of the path $\gamma_{\delta}$ (since the trace is invariant under conjugation),
and so it depends only on the free homotopy class of $\delta$. By a slight abuse of notation we will omit the path $\gamma_{\delta}$ in the
notation and we will write just
\begin{equation*}
 \zeta_{\delta}(\rho,t) = \exp\left(\sum_{m\geq 1} \varepsilon(\delta^m) \Tr\left(\rho\left([\delta^m ]\right)\right) \frac{t^m}{m} \right).
\end{equation*}

\begin{Def}
 The \emph{$\rho$-twisted Lefschetz zeta function of $\phi$} is
 \begin{equation*}
  \zeta_{\rho}(\phi) := \prod_{\delta \mbox{\scriptsize{ simple}}} \zeta_{\delta}\left(\rho,t^{\langle \delta,S \rangle}\right) \in R[[t]].
 \end{equation*}
\end{Def}

The $\rho$-twisted Lefschetz zeta functions enjoy similar properties to those of the commutative ones. In particular a result analogue to Corollary
\ref{Corollary: Zeta functions equals torsion in the Abelian case} holds:
\begin{equation} \label{Equation:  Zeta functions equals torsion in the algebraicallly twisted case}
 \zeta_{\rho}(\phi) = \tau_{\rho}(\To_{\phi})
\end{equation}
where $\tau_{\rho}(\To_{\phi}) = \tau_{\rho}(\To_{\phi},t) $ is the Lin's \emph{$\rho$-twisted Reidemeister torsion} of $\To_{\phi}$ (see \cite{J-W}
and \cite{Lin} for the details). In particular $\zeta_{\rho}(\phi)$ depends only on the conjugacy class of $\rho$ and the isotopy class of $\phi$.
The set $\{\zeta_{\rho}(\phi)\}_{\rho}$ is then a topological invariant of $\To_{\phi}$.

\subsection{Twisted Lefschetz zeta functions associated with bundles} \label{Subsection: Geometric twistings for Lefschetz zeta functions}  $ $ \\
\indent Let $F$ be an oriented compact connected surface (possibly with boundary) and
\begin{equation*}
 F \hooklongrightarrow V \stackrel{\pi}{\twoheadlongrightarrow} \To_{\phi}
\end{equation*}
be a smooth oriented fiber bundle. Let 
\begin{equation} \label{Equation: Geometric monodromy map}
 \mathfrak{m} = \mathfrak{m}_{\pi} :\pi_1(\To_{\phi}) \longrightarrow \mathrm{MCG}(F)
\end{equation}
be the \emph{monodromy homomorphism} of the bundle into the oriented mapping class group of $F$. If $\delta$ is an orbit of
$\phi$ and $\gamma_{\delta}$ is a path from the base point for $\pi_1(\To_{\phi})$ to $\delta$ like in last subsection, then
$\m([\gamma_{\delta} \delta \gamma_{\delta}^{-1}])$ depends only on the
free homotopy class of $\delta$, so that $\m([\delta])$ is well defined. For any simple orbit $\delta$ of $\phi$, let
$$\psi_{\delta}:F \longrightarrow F$$
be a smooth non-degenerate representative of $\m([\delta])$ and, for any $n>1$, set $\psi_{\delta^n} = \psi^n_{\delta}$.
Now, the total space $\pi^*\delta$ of the restriction of $\pi$ to $\delta$ naturally fibers
over $\delta \cong S^1$ with fiber $F$ and we have a diffeomorphism
$$\pi^*\delta \cong \To_{\psi_{\delta}}.$$

\begin{Rmk} \label{Remark: The generator l of the mapping torus projects to delta}
 Given any simple orbit $\delta$ of $\To_{\phi}$, by (\ref{Equation: Homology of a mapping torus}) it follows that
 $$H_1(\pi^*\delta) \cong \mathrm{coker}(\mathbbm{1}-{{\psi_{\delta}}}_1) \oplus \Z_l$$
 where the generator $l$ can be chosen in a way that 
 \begin{equation*}
  \pi_* \circ i_*(l) = [\delta] \in H_1(\To_{\phi})
 \end{equation*}
 where $i:\pi^*\delta \hookrightarrow V$ is the natural inclusion. Moreover
 \begin{equation*}
  \ker(\pi_* \circ i_*|_{H_1(\To_{\psi_{\delta}})})) = \mathrm{coker}(\mathbbm{1}-{{\psi_{\delta}}}_1).
 \end{equation*}
 In particular for any orbit $\gamma$ in $\To_{\psi_{\delta}}$:
 \begin{equation*}
  \pi_* \circ i_*([\gamma]) = \langle \gamma,F\rangle [\delta].
 \end{equation*}
 where $\langle \gamma,F \rangle$ denotes the algebraic intersection number in $\To_{\psi_{\delta}}$ between $\gamma$ and the surface
 $F = F\times \{0\} \subset \To_{\psi_{\delta}}$. Then, for the usual formal variables $t_a$ encoding classes $a\in H_1(\To_{\phi})$,
 there is a natural identification
  \begin{equation*}
  t_{\pi_*\circ i_*([\gamma])} = t_{\langle \gamma,F \rangle[\delta]} = t_{[\delta]}^{\langle \gamma,F \rangle}.
  \end{equation*}
 \begin{comment}
 $\pi_* \circ i_*(l) = \langle \delta,S \rangle \mu$, where
 $i : \pi^*\delta \hookrightarrow V$ is the natural inclusion and $\mu$ is the generator of $H_1(\To_{\phi})$ as in
 (\ref{Equation: Homology of a mapping torus}). Let $t_{\pi_*\circ i_*(l)},\ t_{\mu}$ be formal variables encoding respectively the classes
 $\pi_* \circ i_*(l),\ \mu \in H_1(\To_{\phi})$. Then there is a natural identification
 $$t_{\pi_*\circ i_*(l)} = t_{\mu}^{\langle \delta,S \rangle}.$$
 \end{comment}
\end{Rmk}
%For simplicity we will sometimes consider $\pi^*\delta$ as a submanifold of $V$ avoiding to refer to the natural inclusion map.

\begin{Notation} \label{Notaion:  Homology classes b_V} Similarly to Notation \ref{Notation: Homology classes B_W}, we will often regard $\pi^*\delta$
 as a submanifold of $V$ and
 if $b \in H_1(\pi^*\delta)$, $b_V \in H_1(V)$ will denote the image of $b$ under the homomorphism induced by the inclusion of $\pi^*\delta$ in $V$.
 Moreover we will often assume the identification $\pi^*\delta = \To_{\psi_{\delta}}$ without mentioning it.
\end{Notation}

\begin{Def}  \label{Definition: Local Lefschetz zeta function geometrically twisted}
 Let $(V,\To_{\phi},\pi,F)$ be a surface bundle as above and let $\delta$ be a simple periodic orbit of $\phi$.
 We define the \emph{local $\pi$-twisted Lefschetz zeta function} 
 $$\zeta_{\delta}(\pi,t) \in \Z[[H_1(V)]][[t]]$$
  of $\delta$ by:
 \begin{equation*}
 \zeta_{\delta}(\pi,t) := \prod_{\substack{\gamma \mbox{\scriptsize{ simple}}\\ \mbox{\scriptsize in } \To_{\psi_{\delta}}}} \exp\left(\sum_{m \geq 1} \varepsilon(\delta^{\langle \gamma,F \rangle \cdot m}) \varepsilon(\gamma^m) \frac{\big(z_{[\gamma]_V} \cdot t^{\langle \gamma,F \rangle}\big)^m}{m}\right).
\end{equation*}
 %\begin{equation*}
 % \zeta_{\delta}(\pi) := \prod_{\substack{\gamma \mbox{\scriptsize{ simple}}\\ \mbox{\scriptsize in } \To_{\psi_{\delta}}}} \zeta_{\gamma}\left(z_{[\gamma]_V} \cdot t^{\langle \delta, S \rangle \langle \gamma,F \rangle}\right) \in \Z[[H_1\left(V\right),t]] 
 %\end{equation*}
 % where the formal variable $z_b$ encodes the homology class $b \in H_1(V)$. 
 where the variables $z_d$ keep track of homology classes $d \in H_1(V)$.
\end{Def}

%Remark that $\zeta_{\delta}(\pi,t)$ is a count of orbits with weights only signs and a formal variable encoding 

\begin{Def}  \label{Definition: Lefschetz zeta function geometrically twisted}
 Given $(V,\To_{\phi},\pi,F)$, we define the \emph{$\pi$-twisted Lefschetz zeta function}
 $$\zeta_{\pi}(\phi) \in \Z[[H_1(V)\oplus H_1(\To_{\phi})]][[t]] $$
  of $\phi$ by:
 \begin{equation*}
 \zeta_{\pi}(\phi) := \prod_{\substack{\delta \mbox{\scriptsize{ simple}}\\ \mbox{\scriptsize in } \To_{\phi}}} \zeta_{\delta}\left(\pi,t_{[\delta]}\cdot t^{\langle\delta,S\rangle}\right).
%  \zeta(\pi^*\delta) := \zeta \prod_{\gamma \mbox{\scriptsize{ simple}}}
 \end{equation*}
 where the variables $t_a$ keep track of homology classes $a \in H_1(\To_{\phi})$.
\end{Def}

\subsection{Relation with the non-commutative Lefschetz zeta functions} \label{Subsection: Relation with the algebraic twistings} $ $\\
\indent Given a surface bundle $(V,\To_{\phi},\pi,F)$, we want in some sense the ``richest representation'' of $\pi_1(T_{\phi})$ 
induced by $\m_{\pi}$ into a linear group over some Abelian ring. In what follows we will define this in the fibered case, but the construction can be
carried on in a completely analogous way for a general $3$-manifold (cf. Definition \ref{Definition: Bundled representations} below).

Given $(V,\To_{\phi},\pi,F)$, the associated monodromy map $\m=\m_{\pi}$ in (\ref{Equation: Geometric monodromy map}) induces in homology 
an \emph{algebraic monodromy of index} $i \in \{0,1,2\}$
\begin{equation} \label{Equation: Algebraic monodromy maps}
 \begin{array}{ccccc}
   \m_i = \m_i^{\pi} & : & \pi_1(\To_{\phi}) & \longrightarrow & \mathrm{GL}(H_i(F))\\
                      &   &         \alpha          & \longmapsto     &  \m_i(\alpha)
 \end{array}
\end{equation}
where $r_i = \mathrm{rank}(H_i(F))$ and $\m_i(\alpha)$ is the matrix 
%(with respect to a fixed basis of $H_i(F)$) 
associated to the homomorphism
induced in homology by (any representative of) $\m(\alpha)$.
In order to get all the Abelian information that we can, we
proceed as in Subsection \ref{Subsection: Commutative Lefschetz zeta functions}. Given a class $\alpha \in \pi_1(\To_{\phi})$, we have
$$\pi^*\alpha \stackrel{\mathrm{homeo}}{\cong} \To_{\m(\alpha)}$$
where, with a slight abuse of notation, we keep to use $\alpha$ instead of taking a representative. Let $l_{\alpha}$ denote the generator of
$$H_1(\pi^*\alpha) \cong \mathrm{coker}(\mathbbm{1}-\m_1(\alpha)) \oplus \Z_{l_{\alpha}}$$
as in Remark \ref{Remark: The generator l of the mapping torus projects to delta} and satisfying in particular the relation
$\pi_* \circ i_*(l_\alpha) = [\alpha] \in H_1(\To_{\phi})$.

Then, as we have seen recalling the definition of the total Abelian Lefschetz numbers, for $i\in \{0,1,2\}$, we can go to the universal
Abelian cover of $F$ and consider the \emph{total Abelian algebraic monodromies}
\begin{equation*}
 \begin{array}{ccccc}
   \widetilde{\m}_i & : & \pi_1(\To_{\phi}) & \longrightarrow & \mathrm{GL}\left(r_i,\Z[\mathrm{coker}(\mathbbm{1}-\m_1(\alpha))]\right).\\
        &   &         \alpha          & \longmapsto     &  \widetilde{\m}_i(\alpha)
 \end{array}
\end{equation*}
These induce representations
\begin{equation*}
 \begin{array}{ccccc}
   \widetilde{\m}^V_i  & : & \pi_1(\To_{\phi}) & \longrightarrow & \mathrm{GL}\left(r_i,\Z[H_1(V)]\right).\\
        &   &         \alpha          & \longmapsto     &   z_{({l_{\alpha}}_V)}\cdot i_*(\widetilde{\m}_i(\alpha))
 \end{array}
\end{equation*}
where, as usual, $i_*$ denotes the map induced by the natural inclusion $i:\pi^*\alpha \hookrightarrow V$ and we associate variables $z_{c_V}$
(encoding classes in $H_1(V)$) to classes $c\in H_1(\pi^*\alpha)$. 
Consider now the ring homomorphism
\begin{equation*}
 \begin{array}{cccc}
  p : & \Z[H_1(V)] & \longrightarrow & \Z[H_1(V) \oplus H_1(\To_{\phi})]\\
      &      az_b  & \longmapsto     & a z_b \cdot t_{\pi_*b}
 \end{array}
\end{equation*}
and the induced representations
\begin{equation*}
 \begin{array}{cccc}
  p : & \mathrm{GL}\left(r_i,\Z[H_1(V)]\right) & \longrightarrow & \mathrm{GL}\left(r_i,\Z[H_1(V) \oplus H_1(\To_{\phi})]\right).
 \end{array}
\end{equation*}

%We will denote the $(j,k)$-entry of the matrix
%$\widetilde{\m}^V_i(\alpha)$ by $\widetilde{\m}^V_i(\alpha)_{j,k}$.

%coefficients in H_1(\pi^*\delta) \cong \mathrm{coker}(\mathbbm{1}-{{\psi_{\delta}}_*}_1) \oplus \Z_l,
\begin{Def} \label{Definition: Total Abelian representation induced by pi}
Given the surface bundle $(V,\To_{\phi},\pi,F)$, 
%let $z_{b_{j,k}}$, for $b_{j,k}\in H_1(V)$, be the total $H_1(V)$-component of $\widetilde{\m}^V_i(\alpha)_{j,k}$.
the \emph{total Abelian monodromy representation of $\pi_1(\To_{\phi})$ of index $i \in\{0,1,2\}$ induced by $\pi$} is
\begin{equation*}
  \begin{array}{ccccc}
   \rho^{\pi}_i : & \pi_1(\To_{\phi}) & \longrightarrow & \mathrm{GL}\left(r_i,\Z[H_1(V) \oplus H_1(\To_{\phi})]\right).\\
                  &         \alpha     & \longmapsto     &   p\left(\widetilde{\m}^V_i(\alpha)\right)
 \end{array}
\end{equation*}
\end{Def}

\begin{Thm} \label{Theorem: Relation between geometric and algebraic twistings for the zeta function.}
 Given a surface bundle $(V,\To_{\phi},\pi,F)$ with induced total Abelian monodromy representations $\rho^{\pi}_i$ we have
 $$\zeta_{\pi}(\phi) = \prod_{i=0}^2 (\zeta_{\rho^{\pi}_i}(\phi))^{(-1)^i}$$
 where $\zeta_{\rho^{\pi}_i}(\phi)$ is the $\rho^{\pi}_i$-twisted Lefschetz zeta function of $\phi$.
\end{Thm}

\begin{Cor}
 The $\pi$-twisted Lefschetz zeta function of $\phi$ depends only on the homeomorphism class of $\To_{\phi}$ and the bundle isomorphism
 class of $(V,\To_{\phi},\pi,F)$. In particular $\zeta_{\pi}(\phi)$ does not depend on the choice of the representatives $\psi_{\delta}$
 of $\m([\delta])$.
\end{Cor}

\proof[Proof of Theorem \ref{Theorem: Relation between geometric and algebraic twistings for the zeta function.}] Before going into the computations
we give the basic idea behind the theorem. If $\delta$ is a simple orbit of $\phi$, then the local $\pi$-twisted Lefschetz zeta function
$\zeta_{\delta}(\pi,t)$ of $\delta$ can be interpreted as a kind of ``global'' total Abelian Lefschetz function of the manifold $\To_{\psi_{\delta}}$
for a fixed representative $\psi_{\delta}$ of $\m([\delta])$. Using the definition in (\ref{Equation: Total Abelian Lefschetz numbers}) and the
relations (\ref{Equation: Total Abelian Lefschetz fixed point theorem}) and
(\ref{Equation: Total Abelian Lefschetz zeta function in terms of Lefschetz numbers}) it is possible to express the count of orbits in
$\zeta_{\delta}(\pi,t)$ in terms of traces of powers of the matrices $\widetilde{\m}_i([\delta])$.
%(observe that the algebraic monodromies are well defined on free homotopy class of loops). 

The actual computation requires a bit of attention because we want to keep track of the homology classes of the orbits in $\To_{\psi_{\delta}}$
and also of their projections in $H_1(\To_{\phi})$, which are counted via the more accurate representations $\rho^{\pi}_i$.

Let us begin by translating the orbit count in $\zeta_{\delta}(\pi,t^{\langle \delta,S \rangle})$ in terms of fixed points of the iterates of
the representative $\psi_{\delta}$ of $\m([\delta])$.

\vspace{0.4 cm}

\noindent $\zeta_{\delta}(\pi,t_{[\delta]}\cdot t^{\langle\delta,S\rangle}) =$
\begin{equation*}
 \begin{array}{ll}  
 = & \displaystyle \prod_{\substack{\gamma \mbox{\scriptsize{ simple}}\\ \mbox{\scriptsize in } \To_{\psi_{\delta}}}} \exp\left(\sum_{m \geq 1} \varepsilon(\delta^{\langle \gamma,F \rangle \cdot m}) \varepsilon(\gamma^m) \frac{(z_{[\gamma]_V} )^m \cdot \left(t_{[\delta]}\cdot t^{\langle \delta,S \rangle}\right)^{\langle \gamma,F \rangle \cdot m}}{m}\right)  \\
 = & \displaystyle  \exp\left(\sum_{m \geq 1} \sum_{\substack{\gamma \mbox{\scriptsize{ simple}}\\ \mbox{\scriptsize in } \To_{\psi_{\delta}}}} \varepsilon(\delta^{\langle \gamma,F \rangle \cdot m}) \varepsilon(\gamma^m) z_{[\gamma^m]_V} t_{\pi_*([\gamma^m]_V)} \frac{\left(t^{\langle \delta,S \rangle}\right)^{\langle \gamma,F \rangle\cdot m}}{m}\right) \\
 %= \displaystyle  \exp\left(\sum_{n \geq 1} \sum_{\substack{x \in \mathrm{Fix}(\psi_{\delta}^n)}}  \varepsilon(\delta^n) \varepsilon(x,n) \frac{\left(z_{[\gamma_x]_V}t_{\pi_*([\gamma_x]_V)}t^{\langle \delta,S \rangle}\right)^n}{n}\right) & =\\
 = & \displaystyle  \exp\left(\sum_{n \geq 1}  \sum_{\substack{x \in \mathrm{Fix}(\psi_{\delta}^n)}} \varepsilon(\delta^n) \varepsilon(x,n) \left( z_{[\gamma_x]_V}t_{\pi_*([\gamma_x]_V)}\right)\frac{t^{\langle \delta,S \rangle \cdot n}}{n}\right).
 %= \displaystyle  \exp\left(\sum_{n \geq 1} \varepsilon(\delta^n) \sum_{\substack{x \in \mathrm{Fix}(\psi_{\delta}^n)}}  \varepsilon(x,n) \frac{\left(z_{[\gamma_x]_V}t_{\pi_*([\gamma_x]_V)}t^{\langle \delta,S \rangle}\right)^n}{n}\right) & =\\
 \end{array}
\end{equation*}
Here $n=\langle \gamma,F \rangle \cdot m$ and $\gamma_x$ is the unique $n$-periodic orbit in $\To_{\psi_{\delta}}$ passing through $x$ and has homology class
$[\gamma_x] \in H_1(\To_{\psi_{\delta}})$. By (\ref{Equation: Total Abelian Lefschetz zeta function in terms of Lefschetz numbers}) and 
the definition of $\Le_{\mathrm{A}}$ in (\ref{Equation: Total Abelian Lefschetz numbers}) we get:

\vspace{0.2 cm}

\noindent $ \begin{array}{lll}
  \zeta_{\delta}(\pi,t_{[\delta]}\cdot t^{\langle\delta,S\rangle}) & = & \displaystyle  \exp\left(\sum_{n \geq 1} \varepsilon(\delta^n) \sum_{i=0}^2   (-1)^i \Tr\left(\rho^{\pi}_i([\delta^n])\right)\frac{t^{\langle \delta,S \rangle \cdot n}}{n}\right) \\
     &  = & \displaystyle  \prod_{i=0}^2  \left(\exp\left(\sum_{n \geq 1} \varepsilon(\delta^n) \Tr\left(\rho^{\pi}_i([\delta^n])\right)\frac{t^{\langle \delta,S \rangle \cdot n}}{n}\right)\right)^{(-1)^i} \\
     &  = & \displaystyle  \prod_{i=0}^2 (\zeta_{\delta}(\rho^{\pi}_i,t))^{(-1)^i}.
 \end{array}$

\vspace{0.2 cm}

\noindent The result follows then by multiplying over all the simple periodic orbits $\delta$ of $\phi$.
\endproof

In the proof of last theorem we regarded the contribution to $\zeta_{\pi}(\phi)$ of each simple periodic orbit $\delta$ of $\phi$ as a kind of
total Abelian Lefschetz zeta function of the manifold $\pi^*\delta$, evaluated in the image in $H_1(V) \oplus H_1(\To_{\phi})$ of $H_1(\pi^*\delta)$
via the homomorphism $i_* \oplus (\pi\circ i)_*$, where $i:\pi^*\delta \hookrightarrow V$ is the inclusion.
The only difference with the true total Abelian Lefschetz zeta function $\zeta_{\mathrm{A}}(\pi^*\delta)$ (conveniently evaluated in $H_1(V) \oplus H_1(\To_{\phi})$)
is the sign of the contribution of each orbit $\gamma$ of $\psi_{\delta}$, that in the $\pi$-twisted version is
multiplied by the sign of the power of $\delta$ determined by $\pi_*\circ i_*([\gamma])$. Using the relation in (\ref{Equation: Total Abelian zeta functions and Reidemeister torsions})
we can understand this difference by expressing $\zeta_{\pi}(\phi)$ in terms of the \emph{``local'' Reidemeister torsions} of $(V,\To_{\phi},\pi,F)$
near the orbits of $\phi$. Denote by
$$\tau_V(\pi^*\delta,t) \in \Z[[H_1(V)]][[t]]$$
the series obtained from $\tau_{\mathrm{A}}(\pi^*\delta,t)$ by replacing the variables keeping track of classes $c\in H_1(\pi^*\delta)$ with the
corresponding variables $z_{c_V}$ (using Notation \ref{Notaion:  Homology classes b_V}) encoding classes in $H_1(V)$.

\begin{Prop} \label{Proposition: Geometrically twisted Lefschetz zeta function in terms of local Reidemeister torsions}
  Given a surface bundle $(V,\To_{\phi},\pi,F)$, we have
  \begin{eqnarray*}\zeta_{\pi}(\phi)% & = & \prod_{\substack{\delta \mbox{\scriptsize{ simple}}\\ \mbox{\scriptsize in } \To_{\phi}}} \tau_V\left(\pi^*\delta,\left(\varepsilon(\delta)\varepsilon(\delta^2)\right)\cdot t_{[\delta]}\cdot t^{\langle\delta,S\rangle}\right)^{\varepsilon(\delta^2)}\\                  
                    & = &    \prod_{\substack{\delta \mbox{\scriptsize{ simple}}\\ \mbox{\scriptsize in } \To_{\phi}}}
                       \left\{ \begin{array}{ll}
                        \tau_V(\pi^*\delta,t_{[\delta]}\cdot t^{\langle\delta,S\rangle}) & \mbox{if } \delta \mbox{ elliptic;}\\
                         & \\
                        \dfrac{1}{\tau_V(\pi^*\delta,t_{[\delta]}\cdot t^{\langle\delta,S\rangle})}  & \mbox{if } \delta \mbox{ positive hyperbolic;}\\
                         & \\
                        \dfrac{1}{\tau_V(\pi^*\delta,-t_{[\delta]}\cdot t^{\langle\delta,S\rangle})}  & \mbox{if } \delta \mbox{ negative hyperbolic.}
                      \end{array}\right.
   \end{eqnarray*}
\end{Prop}
\proof
 We prove here only the case of $\delta$ negative hyperbolic, leaving to the reader the analogue computations for the other two cases. For the case
 in hand, by (\ref{Equation: Lefschetz signs of the type of the orbits}) we have $\varepsilon(\delta^m) = (-1)^{m+1}$ and, as in the proof of Theorem
 \ref{Theorem: Relation between geometric and algebraic twistings for the zeta function.}, we get
 $$\zeta_{\delta}(\pi,t_{[\delta]}\cdot t^{\langle\delta,S\rangle}) = \prod_{i=0}^2  \left(\exp\left(\sum_{m \geq 1}  (-1)^{m+1} \Tr\left(\rho^{\pi}_i([\delta^m])\right)\frac{t^{\langle \delta,S \rangle \cdot m}}{m}\right)\right)^{(-1)^i}.$$
 By definition, $\rho^{\pi}_i([\delta^m])$ is the matrix induced on $H_i(\widetilde{F},\Z[H_1(V) \oplus H_1(\To_{\phi})])$
 (via the evaluation induced by $i_*\oplus (\pi\circ i)_*$ on the variables) %for $i:\pi^*\delta \hookrightarrow V$)
 by a representative $\psi_{[\delta^m]}$ of
 $\m([\delta^m])$. Then, using the Taylor series of the formal logarithm $\log$ and the linear algebra relation $\exp(\Tr(\log(\cdot))) = \det(\cdot)$,
 we get:
  \begin{equation*}
   \begin{array}{lll}
    \zeta_{\delta}(\pi,t_{[\delta]}\cdot t^{\langle\delta,S\rangle}) & = & \displaystyle \prod_{i=0}^2  \left(\exp\left(\sum_{m \geq 1} - \Tr\left(\rho^{\pi}_i([\delta^m])\right)\frac{(-t^{\langle \delta,S \rangle})^m}{m}\right)\right)^{(-1)^i} \\
                                                                     & = & \displaystyle \prod_{i=0}^2  \left(\exp\left(\Tr\left(-\sum_{m \geq 1}\frac{(\rho^{\pi}_i([\delta])\cdot (-t^{\langle \delta,S \rangle}))^m}{m}\right)\right)\right)^{(-1)^i} \\
                                                                     & = & \displaystyle \prod_{i=0}^2  \Bigg(\exp\bigg(\Tr\bigg(\log\left(\mathbbm{1}- \rho^{\pi}_i([\delta])\cdot (-t^{\langle \delta,S \rangle})\right)\bigg)\bigg)\Bigg)^{(-1)^i} \\
                                                                     & = & \displaystyle \prod_{i=0}^2  \Bigg(\det\bigg(\mathbbm{1}- \rho^{\pi}_i([\delta])\cdot (-t^{\langle \delta,S \rangle})\bigg)\Bigg)^{(-1)^i} \\
                                                                     & = & \displaystyle \left(\tau_V(\pi^*\delta,-t_{[\delta]}\cdot t^{\langle\delta,S\rangle})\right)^{-1}
  \end{array}
 \end{equation*}
 where the last equivalence follows by applying the relation (\ref{Equation: Total Abelian zeta functions and Reidemeister torsions}) to $\pi^*\delta$
 (with monodromy ${\m}([\delta])$) and replacing as usual the variables encoding the classes $c\in H_1(\pi^*\delta)$ with the variables
 $z_{i_*(c)} \cdot t_{\pi_*\circ i_*(c)}$.
\endproof

Observe that the weight with which we count a simple orbit $\delta$ is, up to eventually taking the reciprocal, just an evaluation of the standard
total Abelian Reidemeister torsion of $\To_{\psi_{\delta}}$.

\begin{Ex}[The trivial example] Consider the disk bundle $(\mathbb{D} \times \To_{\phi},\To_{\phi},\pi,\mathbb{D})$. For any simple orbit $\delta$,
 $\pi^*\delta$ is then a solid torus,
 %(it can be thought of as the complement of the trivial knot in $S^3$),
 whose Reidemeister torsion is
 $$\tau(\mathbb{D}\times S^1,t) = \dfrac{1}{1-t}.$$
 Since here $\pi_*:H_1(V) \rightarrow H_1(\To_{\phi})$ is an isomorphism, we can set all the variables $z_b = 1$ in $\zeta_{\pi}(\phi)$ without loosing information.
 Then Proposition \ref{Proposition: Geometrically twisted Lefschetz zeta function in terms of local Reidemeister torsions} gives
 \begin{equation*}\zeta_{\pi}(\phi) = \prod_{\substack{\delta \mbox{\scriptsize{ simple}}\\ \mbox{\scriptsize in } \To_{\phi}}}
                       \left\{ \begin{array}{ll}
                        \tau(\pi^*\delta,t_{[\delta]}\cdot t^{\langle\delta,S\rangle}) = \dfrac{1}{1-t_{[\delta]}\cdot t^{\langle\delta,S\rangle}} & \mbox{if } \delta \mbox{ ell.;}\\
                         & \\
                        \dfrac{1}{\tau(\pi^*\delta,t_{[\delta]}\cdot t^{\langle\delta,S\rangle})} = 1-t_{[\delta]}\cdot t^{\langle\delta,S\rangle}  & \mbox{if } \delta \mbox{ pos. hyp.;}\\
                         & \\
                        \dfrac{1}{\tau(\pi^*\delta,-t_{[\delta]}\cdot t^{\langle\delta,S\rangle})} = 1+ t_{[\delta]}\cdot t^{\langle\delta,S\rangle}   & \mbox{if } \delta \mbox{ neg. hyp.}
                      \end{array}\right.
  \end{equation*}
  which equals the total Abelian Reidemeister torsion of $\To_{\phi}$ (cf. equations (\ref{Equation: Local Lefschetz Zeta functions of the type of orbits}), 
  (\ref{Equation: Total Abelian Lefschetz zeta function}) and (\ref{Equation: Total Abelian zeta functions and Reidemeister torsions})).
\end{Ex}

\begin{Def} \label{Definition: Bundled representations}
 Given a $3$-manifold $Y$, we say that a representation $\rho: \pi_1(Y) \rightarrow \mathrm{GL}(n,R)$ over an Abelian ring $R$ is a
 (\emph{surface}) \emph{bundled representation} if there exist a smooth surface bundle $(V,Y,\pi,F)$ and a group homomorphism
 $$h : \Z[H_1(V) \oplus H_1(\To_{\phi})] \longrightarrow R$$
 such that $\rho$ is conjugate with $h(\rho_i^{\pi})$
 for some $i\in \{0,1,2\}$, where $\rho_i^{\pi}$ is the total Abelian monodromy representation of $\pi_1(Y)$ of index $i$ induced by $\pi$.
 %and defined  exactly as in Definition \ref{Definition: Total Abelian representation induced by pi}.
\end{Def}

\begin{Ex} \label{Example: All representations on GL 2g Z are bundled} Given a surface bundle $(V,Y,\pi,F)$ with $F$ closed of genus $g$, the
algebraic monodromy maps (defined in (\ref{Equation: Algebraic monodromy maps}))
 $$\m^{\pi}_i  :  \pi_1(Y)  \longrightarrow  \mathrm{GL}(r_i,\Z).$$
 with $r_i = \mathrm{rank}(H_i(F,\Z))$ are clearly bundled representations, since for any $\alpha \in \pi_1(Y)$,
 $\m^{\pi}_i(\alpha) = h(\rho_i^{\pi}(\alpha))$ with
 $$\begin{array}{cccc}
    h  : & \Z[H_1(V) \oplus H_1(Y)] & \longrightarrow & \Z \\
         &     n z_{b}  t_a              & \longmapsto     &  n
   \end{array}$$
 for any $b \in H_1(V)$ and $a \in H_1(Y)$. On the other hand it is a classical result (see for example \cite{F-M}) that the map
 $$(V,Y,\pi,F) \longmapsto \m_{\pi} \in \mathrm{MCG}(F)$$
 that associates to a closed surface bundle the corresponding monodromy homomorphism (\ref{Equation: Geometric monodromy map}) induces a bijection
 between the set of bundle isomorphism classes over $Y$ with fiber $F$ and the set of conjugacy classes of representations of $\pi_1(Y)$ into
 $\mathrm{MCG}(F)$. 
 
 These facts, together with the well known surjection
 $$\begin{array}{ccc}
   \mathrm{MCG}(F) & \twoheadlongrightarrow & \mathrm{GL}(H_1(F,\Z)),\\
       \psi            & \longmapsto     & \psi_*|_{H_1(F,\Z)}
   \end{array}$$
 imply that, for any $Y$ and $g \geq 1$ all representations of  $ \pi_1(Y)$ over $\mathrm{GL}(2g,\Z)$ are bundled.
\end{Ex}

Observe that the total Abelian Reidemeister torsion  of $\To_{\phi}$ as defined in 
(\ref{Equation: Total Abelian zeta functions and Reidemeister torsions}) is the product of three factors. The $i$-th \emph{total Abelian Alexander
polynomial of $\To_{\phi}$} is
$$\Delta_i(\To_{\phi},t) = \det\left(\mathbbm{1} - t(t_{\mu}\widetilde{\phi}_i)\right).$$
Given an orbit $\delta$ of $\phi$, define ${\Delta_i}_V(\pi^*\delta,t) \in \Z[[H_1(V)]][[t]]$ by evaluating in $H_1(V)$
the variables (encoding elements in $H_1(\pi^*\delta)$) of $\Delta_i(\To_{\phi},t)$ as usual, so that
$$\tau_V(\pi^*\delta,t) = \prod_{i=0}^2\left({\Delta_i}_V(\To_{\phi},t)\right)^{(-1)^{i+1}}.$$

\begin{Cor} \label{Corollary: Twisted torsion in terms of local torsions for bundled representations}
 Let $R$ be an Abelian ring and $\rho : \pi_1(\To_{\phi}) \rightarrow \mathrm{GL}(n,R)$ be a bundled representation which is conjugate with
 $h(\rho_i^{\pi})$ for some $i\in\{0,1,2\}$, bundle $(V,\To_{\phi},\pi,F)$ and homomorphism $h$ as in Definition \ref{Definition: Bundled representations}.
 Then:
 \begin{equation*}
  \tau_{\rho}(\To_{\phi},t) = \prod_{\substack{\delta \mbox{\scriptsize{ simple}}\\ \mbox{\scriptsize in } \To_{\phi}}}
                       \left\{ \begin{array}{ll}
                        h\left({\Delta_i}_V(\pi^*\delta,t_{[\delta]}\cdot t^{\langle\delta,S\rangle})\right) & \mbox{if } \delta \mbox{ ell.;}\\
                         & \\
                        \dfrac{1}{h\left({\Delta_i}_V(\pi^*\delta,t_{[\delta]}\cdot t^{\langle\delta,S\rangle})\right)}  & \mbox{if } \delta \mbox{ pos. hyp.;}\\
                         & \\
                        \dfrac{1}{h\left({\Delta_i}_V(\pi^*\delta,-t_{[\delta]}\cdot t^{\langle\delta,S\rangle})\right)}  & \mbox{if } \delta \mbox{ neg. hyp.}
                      \end{array}\right.
 \end{equation*}
 %{\Delta_i}_V(\pi^*\delta,t_{[\delta]}\cdot t^{\langle\delta,S\rangle}) \det\bigg(\mathbbm{1}- \rho^{\pi}_i([\delta])\cdot (-t^{\langle \delta,S \rangle})\bigg)
\end{Cor}

\vspace{0.2 cm}

\begin{Rmk}
 Theorem \ref{Theorem: Relation between geometric and algebraic twistings for the zeta function.}, Proposition
 \ref{Proposition: Geometrically twisted Lefschetz zeta function in terms of local Reidemeister torsions} and Corollary
 \ref{Corollary: Twisted torsion in terms of local torsions for bundled representations} give an answer to Question \ref{Question 2} of the
 introduction, at least for bundled representations of fibered $3$-manifolds. Essentially, for a bundled representation $\rho$ of
 $\pi_1(\To_{\phi})$ there exists a surface bundle over $\To_{\phi}$ for which we can describe the corresponding Lin's non-Abelian Reidemeister
 torsion $\tau_{\rho}(\To_{\phi})$ in terms of the ``local'' Alexander polynomials of the mapping tori (induced by the bundle) above the
 simple orbits of $\phi$. Each of these mapping tori can be considered as a \emph{local model} for the corresponding orbit $\delta$ that encodes
 in its topology informations about the free homotopy class of $\delta$. 
\end{Rmk}

\begin{Rmk} \label{Remark: What are the bundled representations?}
 Here we do not completely study the family of the bundle representations of $\pi_1(Y)$.
 %which are bundled with respect to some surface bundle $(V,Y,\pi,F)$. 
 Example \ref{Example: All representations on GL 2g Z are bundled} (and its analogue for bundles with fibers
 with boundary) gives a relatively big class of examples. We observe that this class can be easily enlarged by using the representations
 $\widetilde{\m}_i^V$ (and not just the weaker $\m_i$) and composed with some representation of $H_1(V)$ that is non-trivial on $\ker(\pi_*|_{H_1(V)})$.
\end{Rmk}

\subsection{Twisted GT and Lefschetz zeta functions} $ $\\
\indent For this subsection, fix a closed symplectic surface $(S,\Omega)$ of genus $g$ and a mapping class $\m \in \mathrm{MCG}(S)$.
Let $\phi$ be a symplectic representative of $\m$ and consider the $4$-manifold 
$$X_{\phi} := S^1 \times \To_{\phi}.$$
$X_{\phi}$ naturally fibers over a torus with fiber $S$ and can be endowed with a symplectic form $\omega$ given by Thurston's Theorem
\ref{Theorem: Thurston}:  
\begin{equation} \label{Equation: Symplectic form on S^1 times mapping torus}
 \omega := ds \wedge dt + \Omega_{X_{\phi}},
\end{equation}
where:\\
\textbullet $\ $ $s$ and $t$ are coordinates for $S^1$ and, respectively, $\frac{[0,+\infty)}{t\sim t+1}$ in $\To_{\phi}$ (with a slight abuse
\phantom{cc} of notation, we identified the coordinates $(s,t)\in S^1 \times S^1$ of the base torus with
\phantom{cc} the corresponding coordinates in $X_{\phi}$);\\
\textbullet $\ $ $\Omega_{X_{\phi}}$ is a closed $2$-form on $X_{\phi}$ that restricts to $\Omega$ over each fiber $F$.\\

We hope that the use of the letter $t$ to denote both the coordinate for $[0,+\infty)$ and the formal variable of Lefschetz zeta functions
and Reidemeister torsions will not induce confusion.

%In this situation it is possible to directly prove the following result (see \cite{I-P 2} or \cite[Section 2.6]{Hu1}).
\begin{Thm}[\cite{I-P 2},\cite{Hu1}] \label{Theorem: GT equals Lefschetz in the total Abelian case} If the genus of $S$ is $g \geq 2$ then
 $$\GT(X_{\phi},\omega) = \zeta_{\mathrm{A}}(\phi)|_{t=1}.$$
\end{Thm}
\proof[Sketch of the proof.]
 Let $J$ be an $\omega$-compatible almost complex structure on $X_{\phi}$ such that $J\partial_s = \partial_t$. Using the condition $g\geq 2$,
 one can show that, for any $A \in H_2(X_{\phi})$, $d_A \geq 0$ if and only if $A$ is a K\"{u}nneth product $A = [S^1] \times a$ for some
 $a \in H_1(\To_{\phi})$, in which case $d_A = 0$. If $\delta$ is any periodic orbit of $\phi$, the surface $S^1 \times \delta$ is $J$-holomorphic
 and these are in fact all the $J$-holomorphic curves in $(X_{\phi},J)$ that we need to consider (cf. \cite[Lemma 2.6]{Hu1}). Clearly the curve
 $$C_{\delta} := S^1 \times \delta$$
 is embedded if and only if $\delta$ is simple. Moreover it is possible to prove that $C_{\delta}$ is regular
 (so that $\epsilon(C_{\delta})$ is well defined) if and only if $\delta$ is non-degenerate, in which case
 $$\epsilon(C_{\delta^n}) = \varepsilon(\delta^n)$$
 for any $n>0$ (see \cite[Lemma 2.7]{Hu1} for details). We have then the following exhaustive correspondence:
 \begin{equation} \label{Equation: Correspondence between types of C_delta and delta}
  \begin{array}{ccl}
   C_{\delta} \mbox{ of type } (+,0) & \longleftrightarrow & \delta \mbox{ elliptic},\\
   C_{\delta} \mbox{ of type } (-,0) & \longleftrightarrow & \delta \mbox{ positive hyperbolic},\\
   C_{\delta} \mbox{ of type } (+,1) & \longleftrightarrow & \delta \mbox{ negative hyperbolic}.
  \end{array}
 \end{equation}
 Comparing (\ref{Equation: Local Lefschetz Zeta functions of the type of orbits}) and (\ref{Equation: Generating functions standard}) normalized
 by (\ref{Equation: Taubes normalization of P_+0}) we have then
 $$P(C_\delta,t) = \zeta_{\delta}(t)$$
 and, with the natural identification $t_{[C_{\delta}]} \equiv t_{[\delta]}$, we obtain:
 $$\GT(X_{\phi},\omega) = \prod_{C_{\delta} \mbox{\scriptsize{ embedded}}} P(C_\delta,t_{[C_{\delta}]}) = \prod_{\delta \mbox{\scriptsize{  simple}}} \zeta_{\delta}(t_{[\delta]}) = \zeta_{\mathrm{A}}(\phi)|_{t=1}.$$
\endproof

The following is essentially Theorem \ref{Theorem: GT twisted = Zeta twisted in intro} in introduction.

\begin{Thm} Let $(V,\To_{\phi},\pi,F)$ be a smooth surface bundle with $F$ closed and let
 $(S^1 \times V,X_{\phi},\Pi,F)$ be the natural bundle induced by the multiplication by $S^1$, where $\Pi := \mathrm{Id} \times \pi$. Endow $X_{\phi}$
 with the symplectic form $\omega$ given in (\ref{Equation: Symplectic form on S^1 times mapping torus}) and suppose that both $F$ and
 $S$ have genus greater than $1$. Then: 
 $$\GT_{\Pi}(X_{\phi},\omega) = \zeta_{\pi}(\phi)|_{t=1}.$$
\end{Thm}
\proof
 By definition of $\GT_{\Pi}(X_{\phi},\omega)$ and the exhaustive correspondence (\ref{Equation: Correspondence between types of C_delta and delta}) we have
 $$\GT_{\Pi}(X_{\phi},\omega) = \prod_{\delta \mbox{\scriptsize{ simple}}} P_{\Pi}(C_\delta).$$
 It is then enough to prove that, for any $\delta$:
 \begin{equation} \label{Equation: P of C_delta equals zeta of delta twisted}
  P_{\Pi}(C_\delta) =  \zeta_{\delta}(\pi,t_{[\delta]})
 \end{equation}
 with the two identifications
 \begin{equation} \label{Equation: Identifications between the variables in the proof of GT_pi = zeta_pi}
  \begin{array}{ll}
  t_{[S^1] \times a} \equiv t_a & \mbox{ for any } a \in H_1(\To_{\phi})\\
  z_{[S^1] \times b} \equiv z_b & \mbox{ for any } b \in H_1(\pi^*\delta).
  \end{array}
 \end{equation}
 between the variables in the definitions \ref{Definition: GT_X for embedded tori} - \ref{Definition: GT_X for covers of embedded tori} and those in definitions
 \ref{Definition: Local Lefschetz zeta function geometrically twisted} - \ref{Definition: Lefschetz zeta function geometrically twisted}. 
 %The proof will easily follow from Proposition \ref{Proposition: Geometrically twisted Lefschetz zeta function in terms of local Reidemeister torsions},
 %Theorem \ref{Theorem: GT equals Lefschetz in the total Abelian case} and the equivalence (\ref{Equation: Total Abelian zeta functions and Reidemeister torsions})
 %between total Abelian Reidemeister torsion and Lefschetz zeta function.
 
 For any $\delta$ simple fix first a symplectic form $\Omega_{\delta}$ on $F$ and a symplectic representative $\psi_{\delta} : F \rightarrow F$ 
 of $\m_{\pi}([\delta]) \in \mathrm{MCG}(F)$. Then we can endow
 $$\Pi^*C_{\delta} = S^1 \times \pi^*\delta \cong X_{\psi_{\delta}}$$
 with a symplectic form $\omega_{\delta}$ like in (\ref{Equation: Symplectic form on S^1 times mapping torus}) that, being of the form
 (\ref{Equation: Symplectic form above the curves}), can be used to compute $P_{\Pi}(C_\delta)$.
 We can then check Equation (\ref{Equation: P of C_delta equals zeta of delta twisted}) case by case using the correspondence in
 (\ref{Equation: Correspondence between types of C_delta and delta}). If $\delta$ is elliptic
 (\ref{Equation: P of C_delta equals zeta of delta twisted}) is satisfied since
 \begin{equation*}
 % \begin{array}{ll}
   P_{\Pi}(C_{\delta}) = \GT_{X_{\phi}}(\Pi^*C_{\delta})  = \tau_V(\pi^*\delta,t_{[\delta]}) = \zeta_{\delta}(\pi,t_{[\delta]}),
 % \end{array}
 \end{equation*}
 where the second equality follows by applying Theorem \ref{Theorem: GT equals Lefschetz in the total Abelian case} to
 $\Pi^*C_{\delta} \cong X_{\psi_{\delta}}$
 with the identifications (\ref{Equation: Identifications between the variables in the proof of GT_pi = zeta_pi}) and the equivalence 
 (\ref{Equation: Total Abelian zeta functions and Reidemeister torsions}),
 %between total Abelian Reidemeister torsion and Lefschetz zeta function
 while the last equality comes from Proposition \ref{Proposition: Geometrically twisted Lefschetz zeta function in terms of local Reidemeister torsions}.
 The proof for $\delta$ positive hyperbolic works in a completely analogous way by just taking the reciprocals in the last equation.
 
 Finally, suppose that  $\delta$ is negative hyperbolic. Then the double cover $\iota$ of $C_{\delta}$ for which
 $\epsilon((C_{\delta})_{\iota}) = -1$ ``doubles $C_{\delta}$ in the $t$-direction of $X_{\phi}$'', so that 
 $$P_{\Pi}(C_{\delta}) = \dfrac{\GT_{X_{\phi}}(\Pi^*C_{\delta})}{\GT_{X_{\phi}}(\Pi^*C_{\delta^2})}.$$
 Observing that $\Pi^*C_{\delta^2} \cong X_{\psi_{\delta}^2}$, reasoning as in the elliptic case  we get
 \begin{equation*}
  P_{\Pi}(C_{\delta}) = \dfrac{\tau_V(\pi^*\delta,t_{[\delta]})}{\tau_V(\pi^*\delta^2,t_{[\delta]}^2)}
 \end{equation*}
 where $\tau_V(\pi^*\delta^2,t)$ is obtained from $\tau_{\mathrm{A}}(\To_{\psi_{\delta}^2},t) = \tau_{\mathrm{A}}(\pi^*\delta^2,t)$ by replacing the variables
 $z_b$, $b\in H_1(\pi^*\delta^2)$, with variables $z_{b_V}$ where $b_V \in H_1(V)$ is
 the image of $b$ via the composition of the homomorphisms $H_1(\pi^*\delta^2) \rightarrow H_1(\pi^*\delta)$ (induced by the double cover projection)
 and $H_1(\pi^*\delta) \rightarrow H_1(V)$ (induced by the inclusion). Then:
 \begin{equation*}
  \begin{array}{lll}
   P_{\Pi}(C_{\delta}) & = & \dfrac{\displaystyle \prod_{i=1}^2 \left(\det\big(\mathbbm{1} - t_{[\delta]}(\widetilde{\m}_i^V)\big)\right)^{(-1)^{i+1}}}{\displaystyle \prod_{i=1}^2 \left(\det\big(\mathbbm{1} - t^2_{[\delta]}(\widetilde{\m}_i^V)^2\big)\right)^{(-1)^{i+1}}}  \\ 
                       &   & \\                     
                       & = & \displaystyle \prod_{i=1}^2 \left( \dfrac{\det\left(\mathbbm{1} - t_{[\delta]}(\widetilde{\m}_i^V)\right)}{\det\left(\mathbbm{1} - t_{[\delta]}(\widetilde{\m}_i^V)\right)\cdot \det\left(\mathbbm{1} + t_{[\delta]}(\widetilde{\m}_i^V)\right)}\right)^{(-1)^{i+1}} \\
                       &   & \\ 
                       & = & \displaystyle \prod_{i=1}^2 \left( \dfrac{1}{\det\left(\mathbbm{1} - (-t_{[\delta]})(\widetilde{\m}_i^V)\right)}\right)^{(-1)^{i+1}} \\
                       &   & \\ 
                       & = & \dfrac{1}{\tau_V(\pi^*\delta,-t_{[\delta]})} = \zeta_{\delta}(\pi,t_{[\delta]}).
 \end{array}
 \end{equation*}
\endproof

The last theorem, together with Theorem \ref{Theorem: Relation between geometric and algebraic twistings for the zeta function.} and Corollary 
\ref{Corollary: Twisted torsion in terms of local torsions for bundled representations}, gives an answer to Question \ref{Question 1} in the
introduction for twisted Reidemeister torsions associated with surface bundles or bundled representations. 
We end the paper with the following two 
natural questions. Given a $4$-dimensional symplectic manifold $(X,\omega)$, fix a smooth surface bundle $(W,X,\pi,F)$ with $F$ closed and $[F] \neq 0$
in $H_2(W,\R)$:
 \begin{enumerate}
  \item Is $\GT_{\pi}(X,\omega)$ independent (eventually up to global shifts on the variables) on the choice of $\omega$?
  \item Observing that $W$ can always be endowed with a symplectic form $\omega_W$ (given by Thurston's theorem), is there
  any relation between $\GT_{\pi}(X,\omega)$ and the higher-dimensional Gromov series of $(W,\omega_W)$? 
 \end{enumerate}

%\begin{QuestionNo} Given $(W,X,\pi,F)$ as in last question, 
 
%\end{QuestionNo}

%\begin{Rmk}
% In this paper we considered only surface bundles because we know that Gromov-Taubes invariants for symplectic $4$-manifolds are smooth topological
% invariants. 
 %From the point of view of higher-dimensional Gromov invariants, 
% It would be interesting to study the possible relations between our twisted Gromov-Taubes invariants associated with surface bundles and the
% higher-dimensional Gromov series of Ionel and Parker (\cite{I-P 2}) for the $6$-dimensional total space of the corresponding bundle (by
% Theorem \ref{Theorem: Thurston} this is always a symplectic manifold).
 
% On the other hand, in the context of Lefschetz zeta functions, we could easily generalize the definitions of the bundle twistings to bundles with
% higher-dimensional fibers, since Reidemeister torsions and Lefschetz zeta functions of higher-dimensional mapping tori are still topological
% invariants. 
%\end{Rmk}

\end{document}